\documentclass[a4paper,10pt,leqno,twoside]{tobart}

\usepackage[english]{babel} \usepackage{inputenc, amsmath, amssymb, latexsym,
  epic, epsfig, rotating, fancyheadings, amsthm, pifont, empheq}

\newcommand{\ldot}{\ensuremath{\textbf{.}}}
\newcommand{\ld}{\ensuremath{,\ldots,}}
\newcommand{\ssq}{\ensuremath{\subseteq}}
\newcommand{\smin}{\ensuremath{\setminus}}
\newcommand{\eps}{\ensuremath{\epsilon}}

\newcommand{\T}{\ensuremath{\mathbb{T}}}

\newcommand{\N}{\ensuremath{\mathbb{N}}} 
\newcommand{\R}{\ensuremath{\mathbb{R}}}
\newcommand{\Z}{\ensuremath{\mathbb{Z}}}
\newcommand{\Q}{\ensuremath{\mathbb{Q}}}

\newcommand{\Leb}{\ensuremath{\mathrm{Leb}}}

\newcommand{\closure}{\ensuremath{\mathrm{cl}}}

\newcommand{\kreis}{\ensuremath{\mathbb{T}^{1}}}
\newcommand{\ntorus}[1][2]{\ensuremath{\mathbb{T}^{#1}}}

\newcommand{\sltr}{\ensuremath{\textrm{SL}(2,\mathbb{R})}}
\newcommand{\torus}{\ensuremath{\mathbb{T}^2}}

\newcommand{\twomatrix}[4]{\ensuremath{\left(\begin{array}{cc} #1 & #2 \\ #3 &
      #4 \end{array}\right)}}

\newcommand{\alphlist}{\begin{list}{(\alph{enumi})}{\usecounter{enumi}\setlength{\parsep}{2pt}
      \setlength{\itemsep}{1pt} \setlength{\topsep}{5pt}
      \setlength{\partopsep}{3pt}}}

\newcommand{\arablist}{\begin{list}{(\arabic{enumi})}{\usecounter{enumi}\setlength{\parsep}{2pt}
          \setlength{\itemsep}{1pt} \setlength{\topsep}{5pt}
          \setlength{\partopsep}{3pt}}}

\newcommand{\romanlist}{\begin{list}{(\roman{enumi})}{\usecounter{enumi}\setlength{\parsep}{2pt}
              \setlength{\itemsep}{1pt} \setlength{\topsep}{5pt}
              \setlength{\partopsep}{3pt}}}

 \newcommand{\listend}{\end{list}}

\newcommand{\bulletlist}{\begin{list}{$\bullet$}{\setlength{\parsep}{2pt}
                \setlength{\itemsep}{1pt} \setlength{\topsep}{5pt}
                \setlength{\partopsep}{3pt}\setlength{\leftmargin}{15pt}}}

\newcommand{\myproof}{\textit{Proof. }}

\newcommand{\foot}{\footnote}

\newcommand{\ncap}{\ensuremath{\bigcap_{n\in\N}}}

\newcommand{\nLim}{\ensuremath{\lim_{n\rightarrow\infty}}}

\newcommand{\ntel}{\ensuremath{\frac{1}{n}}}

\newcommand{\halb}{\ensuremath{\frac{1}{2}}}

\newcommand{\viertel}{\ensuremath{\frac{1}{4}}}

\newcommand{\thx}{\ensuremath{(\theta,x)}}
\newcommand{\thom}{\ensuremath{\theta + \omega}}

\newcommand{\fthx}{\ensuremath{f_{\theta}(x)}}

\newtheoremstyle{tobthm}{3pt}{3pt}{\itshape}{0pt}{\bfseries}{.}{0.5eM}{}
\theoremstyle{tobthm}

\newtheorem{definition}{Definition}[section]
\newtheorem{thm}[definition]{Theorem}

\newtheorem{lem}[definition]{Lemma}
\newtheorem{cor}[definition]{Corollary}  
\newtheorem{prop}[definition]{Proposition}

\newtheorem{conj}[definition]{Conjecture}

\newtheoremstyle{tobrem}{3pt}{3pt}{\normalfont}{0pt}{\bfseries}{.}{0.5em}{}
\theoremstyle{tobrem}

\newtheorem{rem}[definition]{Remark}

\numberwithin{equation}{section}
\numberwithin{figure}{section}

\title{\Large\textsc{Strange nonchaotic attractors in quasiperiodically forced
    circle maps: Diophantine forcing}} \author{T.~J\"ager\thanks{Department of
    Mathematics, TU Dresden, Germany. Email: {\tt
      Tobias.Oertel-Jaeger@tu-dresden.de}}}

\newcommand{\flthx}{\ensuremath{f_{\tau,\theta}(x)}}
\newcommand{\flth}{\ensuremath{f_{\tau,\theta}}}
\newcommand{\nofolge}[1]{\ensuremath{(#1)_{n\in\mathbb{N}_0}}}

\begin{document}

\setlength{\abovedisplayskip}{1.0ex}
\setlength{\abovedisplayshortskip}{0.8ex}

\setlength{\belowdisplayskip}{1.0ex}
\setlength{\belowdisplayshortskip}{0.8ex}

\maketitle 

\abstract{We study parameter families of quasiperiodically forced (qpf) circle
  maps with Diophantine frequency. Under certain ${\cal C}^1$-open
  conditions concerning their geometry, we prove that these families exhibit
  nonuniformly hyperbolic behaviour, often referred to as the existence of
  strange nonchaotic attractors, on parameter sets of positive measure. This
  provides a nonlinear version of results by Young on quasiperiodic
  \sltr-cocycles and complements previous results in this direction which hold
  for sets of frequencies of positive measure, but did not allow for an explicit
  characterisation of these frequencies. As an application, we study a
  qpf version of the Arnold circle map and show that the
  Arnold tongue corresponding to rotation number $1/2$ collapses on an open set
  of parameters.

   The proof requires to perform a parameter exclusion with respect to some twist
   parameter and is based on the multiscale analysis of the dynamics on certain
   dynamically defined critical sets. A crucial ingredient is to obtain good
   control on the parameter-dependence of the critical sets. Apart from the
   presented results, we believe that this step will be important for obtaining
   further information on the behaviour of parameter families like the qpf Arnold
   circle map. }

\section{Introduction}

After the discovery of strange chaotic attractors in two-dimensional dynamical
systems like the H\'enon map \cite{benedicks/carleson:1991}, a natural question
that occurred was to determine the simplest type of smooth
systems that exhibit {\em `strange'} attractors. In particular, it was not
clear whether chaos was a necessary prerequisite for the existence of such 
objects.  Understanding {\em `strange'} in a broad sense as {\em `having a
complicated structure and geometry'} (compare \cite{milnor:1985}), Grebogi {\em
et al} gave a negative answer to this by showing that strange
non-chaotic attractors (SNA) can appear in quasiperiodically forced (qpf)
monotone interval maps \cite{grebogi/ott/pelikan/yorke:1984}. Their argument was
heuristic, but later made rigorous by Keller \cite{keller:1996}.  These findings
prompted further investigations on qpf 1D maps, which have, despite their simple
structure, surprisingly rich dynamics and appear as natural models for physical
systems subject to the influence of two or more external periodic factors with
incommensurate frequencies
\cite{romeiras/etal:1987,ding/grebogi/ott:1989,feudel/kurths/pikovsky:1995}.

For quite a while, studies on the topic were mainly numerical and rigorous
results remained rare. The only exception, apart from the very particular type
of examples in \cite{grebogi/ott/pelikan/yorke:1984,keller:1996}, is the rich
theory of quasiperiodic (qp) \sltr-cocycles and their associated
linear-projective actions. For these systems, the existence of SNA had already
been proved prior to the work of Grebogi {\em et al} by Million\u{s}\u{c}ikov
\cite{millionscikov:1969}, Vinograd \cite{vinograd:1975} and, in a more general
way, Herman \cite{herman:1983}. In this context, the phenomenon is referred to
as the non-uniform hyperbolicity of the cocycle. Due to close relations to the
spectral properties of 1D Schr\"odinger operators with quasiperiodic potential
(see, for example, \cite{avila/krikorian:2004,haro/puig:2006}), there have been
intense efforts to understand the dynamics of qp \sltr-cocycles during the last
three decades (see
\cite{puig:2004,avila/jitomirskaya:2005,avila/jitomirskaya:2010,avila:2010} for
some recent advances). Unfortunately, most methods from this theory cannot
simply be carried over to more general `non-linear' qpf systems, since they
strongly depend on the linear structure and, in many cases, on the close
relations to spectral theory. At the same time, it is also difficult to compare
SNA with the strange attractors appearing in H\'enon-like maps, since on a formal
level these are quite different objects.  Nevertheless, the methods used by
Benedicks and Carleson's in their seminal work on the H\'enon map
\cite{benedicks/carleson:1991} turned out to be equally fruitful for the description
of SNA. Furthermore, the required inductive schemes are easier to implement in
this context, such that one can reasonably hope to elaborate these techniques
further in order to obtain additional insights about the behaviour and dynamics
of parameter families of qpf circle maps. We will come back to this point at the
end of the introduction.

In the context of qpf systems, multiscale analysis and parameter exclusion in
the spirit of Benedicks and Carleson were introduced by Young in
\cite{young:1997}, where she described non-uniformly hyperbolic dynamics in
certain parameter families of qp \sltr-cocycles. The methods were then applied
to qp Schr\"odinger cocycles by Bjerkl\"ov \cite{bjerkloev:2005a}, who also
extended them to show the minimality of the dynamics. These results were so far
restricted to linear-projective systems, but since the original setting in
\cite{benedicks/carleson:1991} is nonlinear it is not too surprising that it was
eventually possible to adapt the techniques to qpf nonlinear models
\cite{jaeger:2009a}.  This allowed to prove the existence of SNA under rather
general conditions.  In \cite{bjerkloev:2005a,jaeger:2009a}, the parameter
exclusion was performed with respect to the forcing frequency. As a result, one
obtains a set of frequencies of positive measure such that the considered system
forced with these frequencies exhibits nonuniformly hyperbolic
dynamics. The drawback is that this does not yield any statement about
a fixed frequency like the golden mean, which is used in most of the numerical
studies on the topic. Our aim here is to close this gap. This is achieved by
performing a parameter exclusion with respect to some other suitable system
parameter. We thus obtain a nonlinear version of the respective results in
\cite{young:1997}, augmented by the minimality of the dynamics. Using a
particular symmetry, we further show that the Arnold tongue corresponding to
rotation number $1/2$ collapses on an open set of parameters. While the collapse
of tongues has already been described in \cite{jaeger:2009a}, the robustness of
this phenomenon seems to be new.

In order to state a qualitative version of our main result, we let ${\cal F} :=
\{f \in \textrm{Diff}^1(\T^2) \mid \pi_1 \circ f =\pi_1\}$, where
$\textrm{Diff}^1(\T^2)$ denotes the group of diffeomorphisms of the two-torus
$\T^2$ and $\pi_1$ is the projection to the first coordinate. Note that for
$F\in{\cal F}$ we have $F\thx = (\theta,f_\theta(x))$ where $f_\theta(\cdot) =
\pi_2\circ F(\theta,\cdot)$, such that we can view $F$ as a collection of {\em
  fibre maps} $(f_\theta)_{\theta\in\kreis}$. Further, we let 
\[ {\cal P} \ = \ \left\{ (F_\tau)_{\tau\in[0,1]} \mid F_\tau\in{\cal F}\
  \forall \tau\in[0,1] \ \textrm{and}\ (\tau,\theta,x)\mapsto F_\tau(\theta,x)
  \textrm{ is } {\cal C}^1\right\} \ 
\]
be the set of differentiable parameter families in ${\cal F}$. The
fibre maps of $F_\tau$ are denoted by $f_{\tau,\theta}$, that is, $F_\tau(\theta,x) =
(\theta,f_{\tau,\theta}(x))$. Finally, we let ${\cal D}(\sigma,\nu)$ be the set
of frequencies $\omega\in\kreis$ that satisfy the Diophantine condition
$d(n\omega,0) > \sigma \cdot |n|^{-\nu} \ \forall n\in\Z\smin\{0\}$.
\begin{thm} \label{t.mr-qualitative} Given any constants $\sigma,\nu >
  0$, there exists a non-empty set ${\cal U} = {\cal U}(\sigma,\nu)
  \ssq {\cal P}$, open with respect to the induced ${\cal
    C}^1$-topology, with the following property:
  
    For all $(F_\tau)_{\tau\in[0,1]} \in {\cal U}$ and all $\omega\in
    {\cal D}(\sigma,\nu)$ there exists a set $\Lambda_\infty(\omega)\ssq [0,1]$
    of positive measure such that for all $\tau\in\Lambda_\infty(\omega)$ the
    qpf circle diffeomorphism
  \[
     f_\tau \ : \ \thx \mapsto (\thom,\flthx)
  \]
  has a unique strange non-chaotic attractor (see Definition~\ref{d.sna}) which
  supports the unique physical measure of the system. Furthermore, the dynamics
  are minimal.
\end{thm}
As in \cite{jaeger:2009a}, we will provide two different quantitative
versions of Theorem~\ref{t.mr-qualitative} which characterise the set
${\cal U}$ in terms of explicit ${\cal C}^1$-estimates.  Since these
conditions are somewhat technical, we postpone the precise statements
to Section~\ref{MainResults} and concentrate on two explicit examples.

The first quantitative result, Theorem~\ref{t.firstversion} below, applies to
the family
\begin{equation} \label{e.projective-action} f_{a,\tau}(\theta,x) \ = \
  \left(\theta+\omega,\frac{1}{\pi}\arctan\left(a^2\tan(\pi
      x)\right)+g(\theta)+\tau\right) \ 
\end{equation}
where $g:\kreis\to\kreis$ is a differentiable function that satisfies some
non-degeneracy condition stated below. For example, one could take
$g(\theta)=\sin(2\pi\theta)$.  If we denote by $R_\phi$ the rotation matrix with
angle $2\pi\phi$, then $f_{a,\tau}$ is the is the projective action of the qp
\sltr-cocycle given by
\begin{equation}
  \label{e.cocycle-example}
  A(\theta) \ = \ R_{g(\theta)+\tau} \cdot \twomatrix{a}{0}{0}{1/a} \ .
\end{equation}
For this particular system, we obtain the following statement.
\begin{cor}[to Theorem~\ref{t.firstversion} below]
  \label{c.cocycle-example} Suppose $g:\kreis\to\kreis$ is a differentiable
  function and there exists a finite set $\Omega_0$ such that for all
  $\tau\in\kreis\smin \Omega_0$ the set $Z(\tau)=\{\theta\in\kreis\mid
  g(\theta)+\tau=\halb \}$ is finite and $g'$ takes distinct and non-zero values
  at different points of $Z(\tau)$.
  
  Then for all $\sigma,\nu>0$ there exists
  $a_*=a_*(g,\sigma,\nu)>0$ with the following property: for
  all $\omega\in {\cal D}(\sigma,\nu)$ and all $a\geq a_*$
  there exists a set $\Lambda_\infty(a,\omega) \ssq \kreis$ of
  positive measure such that for all
  $\tau\in\Lambda_\infty(a,\omega)$ the map $f_{a,\tau}$
  given by (\ref{e.projective-action}) has a unique SNA and minimal
  dynamics. Further $\Leb_{\kreis}(\Lambda_\infty(a,\omega))$
  goes to $1$ as $a \to \infty$. 

  The same result applies to any sufficiently small ${\cal
    C}^1$-perturbation of the parameter family
  (\ref{e.cocycle-example}).
\end{cor}
This statement follows from Theorem~\ref{t.firstversion} by some standard
estimates. Since our main focus lies on the qpf Arnold circle map, we refer the
reader to \cite[Section 3.8]{jaeger:2009a} for details. We also note that the existence of an SNA for
(\ref{e.projective-action}) is equivalent to the non-uniform hyperbolicity of
the cocycle (\ref{e.cocycle-example}) \cite{haro/puig:2006,jaeger:2006a}. Hence,
the result can be viewed as a perturbation-persistent version of \cite[Theorem
2]{young:1997}.  

The second quantitative version of Theorem~\ref{t.mr-qualitative}, stated as
Theorem~\ref{t.mr-quantitative2} below, is tailor-made for the application to
the qpf Arnold circle map
\begin{equation}
  \label{e.arnold}
  f_{a,b,\tau}(\theta,x) \ = \ \left(\theta+\omega,x+\tau+\frac{a}{2\pi}
  \sin(2\pi x) + g_b(\theta)\right) 
\end{equation}
with forcing function $g_b$ depending on some additional parameter $b$. The
geometry of (\ref{e.arnold}) is quite different to that of the previous example,
since unlike in (\ref{e.projective-action}) the hyperbolicity on the single
fibres is limited (the slope of the fibre maps $f_{a,b,\tau,\theta}$ remains
bounded by $2$ in the invertible regime $|a|\leq 1$). In order to make up for
this, the forcing function $g_b$ must have a particular shape that can be pushed
to some extreme by adjusting the parameter $b$. General conditions for the
family $g_b$ can be deduced from Theorem~\ref{t.mr-quantitative2}. (See also
Remark~\ref{r.forcing-structure}.) Here, we concentrate again on an explicit
example.
\begin{cor}[to Theorem~\ref{t.mr-quantitative2} below]
  \label{c.arnold} 
  Let 
   \begin{equation}
  \label{e.ff1}
  g_b(\theta) \ = \ \arctan(b\sin(2\pi\theta))/\pi 
  \qquad (b\in\R) \   . 
\end{equation}
Then for all $\sigma,\nu>0$ and all $a>0$ there exists
$b_*=b_*(\sigma,\nu,a)>0$ with the following property:

  For all $\omega\in {\cal D}(\sigma,\nu)$ and all $b\geq b_*$
  there exists a set $\Lambda_\infty(a,b,\omega) \ssq \kreis$ of
  positive measure such that for all
  $\tau\in\Lambda_\infty(a,b,\omega)$ the map
  $f_{a,b,\tau}$ given by (\ref{e.arnold}) has a unique SNA
  and minimal dynamics. 
\end{cor}
Apart from the restrictions on the forcing function coming from the lack of
hyperbolicity, a further reason for the particular choice of $g_b$ in
(\ref{e.ff1}) is a special symmetry which appears at $\tau=\halb$.  On the one
hand, the lift $F$ of the map $f_{a,b,\halb}$ satisfies the relation
\begin{equation}
  \label{e.forcing-symmetry}
  F_{\theta}(-x) \ = \ 1-F_{\theta+1/2}(x)  \ ,
\end{equation}
and it can be easily seen that this forces the rotation number\foot{See
  Section~\ref{RotNum} for the definition of the fibred rotation number of a qpf
  circle homeomorphism.}  $\rho(f_{a,b,\halb})$ to be exactly $\halb$. On the
other hand, the map $g_b+\halb$ takes values close to \halb\ only on two
intervals $I_0$ and $I_0+\halb$ around $0$ and \halb, respectively.  These two
intervals play a special role in the multiscale analysis, since they define the
critical sets on the first level of the inductive scheme. Furthermore, as a
consequence of (\ref{e.forcing-symmetry}) the fact that the $n$-th critical
region consists of exactly two intervals $I_n$ and $I_n+\halb$ will
remain true on all levels of the induction. This allows to control the return
times of the critical regions directly by using only the Diophantine condition,
and no parameters have to be excluded in order to avoid fast returns. In other
words, in this particular situation the multiscale analysis can be performed
without any parameter exclusion.  As a consequence, we obtain the following.
\begin{cor}[to Theorem~\ref{t.mr-quantitative2} below]
  \label{t.arnold-half}
  Suppose $g_b$ is chosen as in (\ref{e.ff1}). Then for all $\sigma,\nu >
  0$ and $a>0$ there exists $b_*=b_*(\sigma,\nu,a)$ such that
  for all $\omega\in{\cal D}(\sigma,\nu)$ and $b>b_*$ the map
  $f_{a,b,\halb}$ has a unique SNA and minimal dynamics.
  $b_*(\sigma,\nu,\ldot)$ can be chosen constant on compact subsets of
  $(0,1)$.
\end{cor}
This result has further consequences for the structure of the {\em Arnold tongues}
\begin{equation} \label{e.tongues}
 A_\rho \ = \ \left\{(a,b,\tau) \in [0,1]\times\R\times[0,1] \mid
\rho(f_{a,b,\tau})=\rho\right\} \
\end{equation}
and the associated {\em mode-locking plateaus}
\begin{equation}
  \label{e.ml-plateaus}
  P_{a,b,\rho} \ = \ \{\tau \in [0,1] \mid \rho(f_{a,b,\tau})=\rho\} \ ,
\end{equation}
where $\rho(f_{a,b,\tau})$ denotes the fibred rotation number of
$f_{a,b,\tau}$. We say a mode-locking plateau $P_{a,b,\rho}$ is {\em collapsed}
if it consists of a single point. It is known that $P_{a,b,\rho}$ is collapsed
for all $\rho\notin\Q+\Q\omega$ \cite{bjerkloev/jaeger:2009}, and we
implicitly assume that $\rho$ belongs to the module $\Q+\Q\omega$ whenever we
speak of collapsed or non-collapsed plateaus. A tongue $A_\rho$ is said to be
{\em collapsed} at $(a,b)$ if $P_{a,b,\rho}$ is collapsed.  Minimal dynamics
imply the collapse of a tongue, in the sense that whenever $f_{a,b,\tau}$ is
minimal the tongue $A_\rho$ with $\rho=\rho(f_{a,b,\tau})$ is collapsed at
$(a,b)$ (see Proposition~\ref{p.mode-locking}).  Hence, the tongue corresponding to
rotation number \halb\ is collapsed for all the parameters satisfying the
assertions of Corollary~\ref{t.arnold-half}.
\begin{cor} \label{c.tongue-collapse} Suppose $g_b$ is chosen as in
  (\ref{e.ff1}). Then for all $\sigma,\nu > 0$ there exists an open set
  $B(\sigma,\nu) \ssq (0,1)\times\R^+$ such that for $\omega\in D(\sigma,\nu)$
  and forcing function $g_b$ as in (\ref{e.ff1}) the Arnold tongue $A_\halb$ is
  collapsed at all $(a,b)\in B(\sigma,\nu)$.
\end{cor}
In \cite{jaeger:2009a}, it was shown in a similar way that $A_0$
collapses on sets of parameters $(a,b)$ of positive measure,
and the methods employed there would yield the same result for
$A_{\halb}$.  Hence, the new point here is the robustness of
this phenomenon, that is, the openness of the set $B$ in
Corollary~\ref{c.tongue-collapse}.

As mentioned above, there are many further open problems concerning the
behaviour of parameter families like (\ref{e.projective-action}) or
(\ref{e.arnold}). Probably the most prominent one is the question whether the
rotation number as a function of the twist parameter is a `devils staircase',
meaning that the union of non-collapsed mode-locking plateaus is dense in the
parameter interval. This is true for the unforced Arnold circle map. For qpf
systems, existing results are again restricted to qp \sltr-cocycles. A
particular case is the projective action of the Schr\"odinger cocycle associated
to the so-called almost-Mathieu operator, for which the question became known as
the `Ten Martini Problem'. Recently it has been answered positively in full
generality, meaning for all parameters and all irrational forcing frequencies,
by Avila and Jitomirskaya \cite{avila/jitomirskaya:2005} (after previous
contributions by B\'ellisard and Simon \cite{bellissard/simon:1982} and Puig
\cite{puig:2004}). For the qpf Arnold circle map, still no rigorous results
exist. Moreover the numerical findings are ambiguous. On the one hand, a devils
staircase has been reported for some parameters regions in
\cite{ding/grebogi/ott:1989}. On the other hand the authors of
\cite{stark/feudel/glendinning/pikovsky:2002} numerically detect parameters for
which the 0-tongue is collapsed (a fact which is backed up by rigorous results
in \cite{jaeger:2009a} and, replacing $0$ by $\halb$, also by
Corollary~\ref{c.tongue-collapse}) and report that for these parameters the
mode-locking plateaus vanish and the rotation number strictly increases over a
whole interval. In contrast to this, we believe that a further elaboration of
the presented techniques should allow to prove the following.
\begin{conj} \label{conjecture}
  The set $\Lambda_\infty$ in Corollary~\ref{c.arnold} can be chosen such that
  it is contained in the closure of the union of non-collapsed mode-locking
  plateaus. The same is true for the parameter $\tau=\halb$ in the situation of
  Corollary~\ref{c.arnold}.
\end{conj}
In fact, what should be true is that all parameters $\tau$ for which the
`slow-recurrence conditions' \ref{e.X'n} and \ref{e.Y'n} introduced in the
multiscale analysis scheme below are satisfied can be approximated by
non-collapsed mode-locking plateaus. While this would not formally disprove the
conjecture made in \cite{stark/feudel/glendinning/pikovsky:2002} (since our
methods do not apply to the forcing function considered there), it would provide
strong evidence for the fact that the observation is a numerical
artifact. Furthermore, it could be a first step towards proving the existence of
a devils staircase. Apart from the
intrinsic interest of the above results, the hope to make further progress in
this direction is one of the main motivations for the presented work. 

Concerning the proofs, we will be able to rely to a great extent on the previous
construction in \cite{jaeger:2009a}. In particular, the core part of the proof,
which is the multiscale analysis for the dynamics of a fixed map under some
non-recurrence conditions on certain dynamically defined critical sets, remains
valid and can be used for our purpose without any modifications. We will
therefore be able to concentrate almost exclusively on those aspects of the
proof which differ from the previous one. The only drawback of this is that the
present paper is not self-contained, but depends on a number of statements and
estimates in \cite{jaeger:2009a}. However, as redoing all arguments would only
result in an undue length of the paper and render the decisive differences in
comparison to the previous construction less visible, this seems to be the
appropriate way to proceed. In order not to leave the reader without any
guidance, we will briefly motivate the used statements on a heuristic level.
\smallskip

\noindent {\bf Acknowledgements.} I would like to thank Hakan Eliasson for
inspiring discussions on the subject. This work was carried out in the
Emmy-Noether-Group 'Low-dimensional and non-autonomous Dynamics', which is
supported by the Grant Ja 1721/2-1 of the German Research Council (DFG).

\section{Notation and Preliminaries} \label{NotationPreliminaries}

\subsection{Notation.}\label{Notation}
Given $a,b\in\kreis$, we denote by $[a,b]$ the positively oriented arc from $a$
to $b$. The same notation is used for open and half-open intervals. We write
$b-a$ for the length of $[a,b]$, whereas the Euclidean distance between $a$ and
$b$ will be denoted by $d(a,b)$.  The derivative with respect to a variable
$\xi$ will be denoted by $\partial_\xi$.  On any product space, $\pi_i$ will
denote the projection to the $i$-th coordinate. Quotient maps like the canonical
projections $\R\to\T^1=\R/\Z$, $\R^2\to\T^2=\R^2/\Z^2$ or $\T^1\times\R\to\T^2$
will all be denoted by $\pi$.

If $I(\tau) = (a(\tau),b(\tau))$ is an interval that depends on some parameter
$\tau\in\R$, then we say $I$ is differentiable in $\tau$ if this is true for
both endpoints $a$ and $b$. In this case we write
\[
|\partial_\tau I(\tau)| \ = \ \max\{|\partial_\tau
a(\tau)|,|\partial_\tau b(\tau)|\} \ .
\]
If $I^\iota(\tau)=(a^\iota(\tau),b^\iota(\tau))$ and
$I^\kappa(\tau)=(a^\kappa(\tau),b^\kappa(\tau))$ are two disjoint
intervals depending both on $\tau$, then we write
\[
D_\tau(I^\iota(\tau),I^\kappa(\tau)) \ > \ \eta 
\]
if there holds $\partial_\tau(y(\tau)-x(\tau)) > \eta$ for all possible
choices $x(\tau) = a^\iota(\tau),b^\iota(\tau)$ and $y(\tau) =
a^\kappa(\tau),b^\kappa(\tau)$. We write
\[
|D_\tau(I^\iota(\tau),I^\kappa(\tau))| \ > \ \eta
\]
if either $D_\tau(I^\iota(\tau),I^\kappa(\tau)) > \eta$ or
$D_\tau(I^\kappa(\tau),I^\iota(\tau))>\eta$. In other words,
$|D_\tau(I^\iota(\tau),I^\kappa(\tau))| > \eta$ means that the two intervals
move with speed $>\eta$ relative to each other. 

\subsection{SNA in qpf systems.}
\label{SNA_Preliminaries}

We say a continuous map $f:\torus\to\torus$ is a {\em qpf circle
homeomorphism} if it has skew product structure of the form
\begin{equation}
  \label{e.qpf-circlediff}
  f(\theta,x) \ = \ (\theta+\omega,f_\theta(x)) 
\end{equation}
with irrational $\omega\in\T^1$. The maps $f_\theta(x)=\pi_2\circ
f(\theta,x)$ are called {\em fibre maps} and we write
$f^n_\theta(x)=\pi_2\circ f^n(\theta,x)$ for the fibre maps of the
iterates. An invariant graph of $f$ is a measurable function
$\varphi:\kreis\to\T^1$ that
satisfies
\begin{equation}
  \label{e.inv-graph}
  f_\theta(\varphi(\theta)) 
  \ = \ \varphi(\theta+\omega) \ .
\end{equation}
The corresponding point set $\Phi=\{(\theta,\varphi(\theta))\mid
\theta\in\kreis\}$ will equally be called an invariant graph. We note that in
general multi-valued invariant graphs have to be taken into account as
well. However, since in the situation we consider only single-valued invariant
graphs occur, we restrict to this simple case. (The general definitions can be
found in \cite{jaeger:2009a}.) 

To any invariant graph, an $f$-invariant ergodic measure $\mu_\varphi$ can be
assigned by
\begin{equation}
  \label{e.graphmeasure}
  \mu_\varphi(A) \ = \ \Leb_{\kreis}\left(\pi_1(A\cap \Phi)\right) \ .
\end{equation}
If all fibre maps are ${\cal C}^1$ and the derivative $f_\theta'(x)$ is strictly
positive and depends continuously on $\thx$, we speak of a qpf circle
diffeomorphism.  In this case, the {\em (vertical) Lyapunov exponent} of an
invariant graph is defined as
\begin{equation}
  \label{e.LE}
  \lambda(\varphi) \ = \ \int_{\kreis} \log | 
  \partial_x f_\theta(\varphi(\theta))| \ d\theta \ ,
\end{equation}
In the particular context of qpf systems, SNA are now defined as
follows.
\begin{definition}
  \label{d.sna} A non-continuous invariant graph with negative
  Lyapunov exponent is called a {\em strange nonchaotic attractor
    (SNA)}. A non-continuous invariant graph with positive Lyapunov
  exponent is called a {\em strange nonchaotic repeller (SNR)}.
\end{definition}
\begin{rem}
  It should be said at this point that it it difficult to match this very
  specific definition of SNA with a general concept of strange attractors, as
  discussed for example in \cite{milnor:1985}. For instance, an attractor is usually
  understood to be a compact invariant set, but the point set associated to an
  SNA in the above sense is non-compact due to the non-continuity of the
  invariant graph. One could consider the closure of this set instead, but in
  the situations we describe this will be the whole two-torus, which cannot
  reasonably be called a `strange' object.  However, although the terminology
  might therefore be considered somewhat unfortunate, it has already been used
  for almost three decades in most of the vast physics literature on the
  topic. We therefore prefer to keep with it, simply taking it as a technical
  term specific to the theory qpf systems.

  We also note that due to the negative Lyapunov exponent an SNA attracts a
  positive measure set of initial conditions and therefore carries a physical
  measure.
\end{rem}
A convenient criterion for the existence of SNA involves pointwise Lyapunov
exponents, forwards and backwards in time. These are given by
\begin{equation} \label{e.pointwiseLE}
\lambda^\pm(\theta,x) \ = \ \limsup_{n \to \infty} \ntel |\log \partial_x
f^{\pm n}_\theta(x)|\ .
\end{equation}
The orbit of a point $(\theta,x) \in \ntorus$ with
$\lambda^\pm(\theta,x)>0$ is called a \emph{sink-source-orbit}. The
existence of such orbits implies the existence of SNA. 
\begin{prop}[\cite{jaeger:2009a}] \label{prop:sinksourcesna} Suppose
  $f$ is a quasiperiodically forced circle diffeomorphism which has a
  sink-source-orbit. Then $f$ has both a SNA and a SNR.
\end{prop}
In the particular case of the Harper map, the
existence of a sink-source-orbit is equivalent to Anderson
localisation for the corresponding almost-Mathieu operator (see, for
example, \cite{haro/puig:2006} or \cite[Section 1.3]{jaeger:2006a}).

\subsection{The fibred rotation number and mode-locking.} \label{RotNum} If a
qpf circle homeomorphism $f$ is homotopic to the identity on \torus, it has a
continuous lift $F:\kreis\times\R \to \kreis \times \R$ of the form
$F(\theta,x) = (\theta+\omega,F_\theta(x))$. In this case, the limit
\begin{equation}
  \label{e.rotnum-def}
  \rho(F) =  \nLim (F^n_\theta(x)-x)/n
\end{equation}
exists and is independent of $\thx$ \cite[Section 5.3]{herman:1983}. $\rho(f) :=
\rho(F) \bmod 1$ is called the {\em (fibred) rotation number} of $f$. If the
rotation number remains constant under all sufficiently small ${\cal
  C}^0$-perturbations, we speak of {\em mode-locking}. The mechanism for
mode-locking has been clarified in \cite{bjerkloev/jaeger:2009}. For our
purposes we need the following two consequences.
\begin{prop}[\cite{bjerkloev/jaeger:2009,jaeger:2009a}]
  \label{p.mode-locking}
  Suppose $f$ is a qpf circle homeomorphism and let $f_\eps=R_\eps\circ f$, where
  $R_\eps(\theta,x) = (\theta,x+\eps)$. Then the following holds. \alphlist
\item If $\rho(f) \notin \Q+\Q\omega$ then $\eps\mapsto \rho(f_\eps)$ is
  strictly increasing in $\eps=0$. 
\item If $f$ is minimal then $\eps\mapsto \rho(f_\eps)$ is strictly increasing
  in $\eps=0$. \listend
\end{prop}
In other words, mode-locking cannot occur in situations (a) and (b).

\section{Statement of the main results} \label{MainResults}

The explicit ${\cal C}^1$-open conditions characterising the set ${\cal U}$ in
Theorem~\ref{t.mr-qualitative} are not too complicated if each one is considered
by itself, but altogether they form a rather long list. We therefore prefer not
to include them in Theorem~\ref{t.firstversion}, but to state them separately
before.  We also note that conditions
(\ref{eq:Cinvariance})--(\ref{eq:crossing}) below are precisely those used in
\cite{jaeger:2009a}, whereas the Diophantine condition (\ref{eq:Diophantine})
and the assumptions (\ref{eq:bounddlambda})--(\ref{e.d}) on the dependence on
the twist parameter are new. Let $\Lambda\ssq[0,1]$ be an open interval.

\smallskip

\noindent {\em I. Diophantine condition.} First, recall that $\omega\in{\cal
  D}(\sigma,\nu)$ just means that
\begin{equation}
  \label{eq:Diophantine} \tag{${\cal A}0$}
  d(n\omega,0) \ > \ \sigma\cdot |n|^{-\nu} \quad \forall n\in\Z\smin\{0\} \ .
\end{equation}

\noindent {\em II. Critical regions.} Let $E=[e^-,e^+]$ and $C=[c^-,c^+]$ be two
non-empty, compact and disjoint subintervals of \kreis. We will assume that for
all $\tau\in\Lambda$ there exists a finite union ${\cal I}_0(\tau) \ssq
\kreis$ of ${\cal N}$ disjoint open intervals $I_0^1(\tau) \ld I_0^{\cal
  N}(\tau)$ (the {\em `critical regions'}) such that
\begin{equation}  \label{eq:Cinvariance} \tag{${\cal A}1$}
  \flth(\mbox{cl}(\kreis \smin E)) \ \ssq \ \mbox{int}(C) \ \ \ \ \ \forall
  \theta \notin {\cal I}_0(\tau) \ .
\end{equation}
Note that this implies
\begin{equation} \label{eq:Einvariance} \tag{${\cal A}1'$}
  \flth^{-1}(\mbox{cl}(\kreis \smin C)) \ \ssq \ \mbox{int}(E) \ \ \ \ \ \forall
  \theta \notin {\cal I}_0(\tau)+\omega \ .
\end{equation} 

\noindent \emph{III. Bounds on the derivatives.} Concerning the derivatives of
the fibre maps $\flth$, we will assume that for given $\alpha > 1$ and $p\in\N$
we have
\begin{equation} \tag{${\cal A}2$}
  \label{eq:bounds1}  \alpha^{-p} \ < \ \partial_x\flthx \ < \ \alpha^p \hspace{2eM} \quad
  \forall \thx \in \ntorus \ ;
\end{equation}
\begin{equation} \label{eq:bounds2} \hspace{3eM}  \partial_x\flthx \ > \ \alpha^{2/p}
  \hspace{4.1eM} \quad \forall \thx \in \kreis \times  \tag{${\cal A}3$}
  E \ ; \end{equation}
\begin{equation} \label{eq:bounds3} \hspace{3eM} \partial_x\flthx \ < \
  \alpha^{-2/p} \hspace{4.1eM} \quad \forall \thx \in \kreis \times C \ .
  \tag{${\cal A}4$}
 \end{equation}
 Further, we fix $S>0$ such that
\begin{equation} \label{eq:bounddth}  \tag{${\cal A}5$}
   |\partial_\theta \flthx| \ < \ S \ \ \ \ \ \forall \thx \in
  \ntorus \ .
\end{equation}

\noindent \emph{IV. Transversal Intersections.} The significance of the critical
region ${\cal I}_0(\tau)$ is the fact that due to (\ref{eq:Cinvariance}) this is
the only place where the attracting and the repelling region can `mix up'. In
order to ensure that the intersections of $f_\tau({\cal I}_0(\tau) \times C)$
and $({\cal I}_0(\tau)+\omega)\times E$ are `nice' (transversal in an
appropriate sense), we will assume that
\begin{equation} \label{eq:crossing}  \tag{${\cal A}6$}
  \begin{array}{l} \exists!\theta_\iota^1 \in I_0^\iota(\tau) \textrm{ with  }
    f_{\tau,\theta_\iota^1}(c^+) = e^- \textrm{ \ and } \\ \exists! \theta_\iota^2 \in
    I_0^\iota(\tau) \textrm{ with } f_{\tau,\theta_\iota^2}(c^-) = e^+ \ . \end{array}
\end{equation}
This ensures that the image of $I_0^\iota(\tau) \times C$ crosses the strip
$(I^\iota_0(\tau)+\omega) \times E$ exactly once and not several times. In order
to control the slope of these strips, we will assume that
\begin{equation} 
  |\partial_\theta \flthx| \ > \ s \ \ \ \ \ \forall \thx \in {\cal I}_0(\tau)
  \times \kreis \ \label{eq:s} \tag{${\cal A}7$}
\end{equation}
for some constant $s$ with $0 < s < S$. \smallskip

\noindent {\em V. Dependence on $\tau$.}
First, we fix an upper bound $L$ on $|\partial_\tau
\flthx|$, that is,
\begin{equation} \label{eq:bounddlambda} \tag{${\cal A}8$}
  |\partial_\tau \flthx| \ < \ L\ \ \ \ \ \forall \thx \in \ntorus \ .
\end{equation}
Secondly, we assume that all connected components $I^\iota_0(\tau) =
(a^\iota_0(\tau),b^\iota_0(\tau))$ of ${\cal I}_0(\tau)$ are differentiable with
respect to $\tau$ and that for some constant $\eta>0$ we have
\begin{equation}
  \label{e.I_0-derivative} \tag{${\cal A}9$}
  |D_\tau(I^\iota_0(\tau),I^\kappa_0(\tau))| \ > \ \eta/2 \quad 
  \forall \iota \neq \kappa \in [1,\cal N] \ . 
\end{equation}
Finally we will need an assumption which ensures that condition
(\ref{e.I_0-derivative}) holds in a similar way for the higher order critical
regions ${\cal I}^\iota_n$ defined later on. This is actually the crucial point
of the construction, which allows to adapt the parameter exclusion scheme from
\cite{bjerkloev:2005a,jaeger:2009a} to the considered problem.  For any $\iota
\in [1,{\cal N}]$ we let
\begin{equation}\nonumber
  Q^\iota(\tau) \ := \ \left\{ \partial_\tau \flthx / \partial_\theta \flthx 
    \left|  \ \thx \in I^\iota_0 \times \kreis \right. \right\} \ .
\end{equation}
Then we assume that there holds
\begin{equation} \label{e.d} \tag{${\cal A}10$}
  d(Q^\iota(\tau),Q^\kappa(\tau)) \ > \ \eta \quad \forall \iota \neq \kappa
  \in [1,\cal N] \ .
\end{equation}

\begin{thm} \label{t.firstversion} Let $\omega \in {\cal
    D}(\sigma,\nu)$ and $\delta > 0$. Suppose that $\Lambda \ssq
  [0,1]$ is an open interval and $(F_\tau)_{\tau\in[0,1]} \in {\cal
    P}$ is such that conditions (\ref{eq:Cinvariance})--(\ref{e.d})
  are satisfied on $\Lambda$. Let $\eps_0 := \sup_{\iota\in[1,{\cal
      N}],\tau\in\Lambda}|I^\iota_0(\tau)|$. Then there exist constants
  $\alpha_* = \alpha_*(\sigma,\nu,{\cal N},p,S,s,L,\eta,\delta)$ and $\eps_* =
  \eps_*(\sigma,\nu,{\cal N},p,S,s,L,\eta,\delta)$ such that the following
  holds.

  If $\alpha > \alpha_*$ and $\eps_0 < \eps_*$, then there exists a set
  $\Lambda_\infty=\Lambda_\infty(\omega)\ssq\Lambda$ of measure
  $\Leb(\Lambda_\infty) > \Leb(\Lambda) - \delta$ such that for all $\tau \in
  \Lambda_\infty$ the qpf circle diffeomorphism
  \begin{equation}
    f_\tau \ : \ \thx \mapsto (\thom,\flthx)
  \end{equation}
  satisfies
 \begin{equation} \tag*{$(*)$}
    \label{e.star}
    \left\{
      \begin{array}{cl}
        (*1) & f_\tau \textrm{ has a SNA } 
        \varphi^-  \textrm{ and a SNR } \varphi^+;\vspace{1ex} \\ \vspace{1ex}
        (*2) & \varphi^- \textrm{ and } \varphi^+ \textrm{ are one-valued and the only
          invariant graphs of } f;\\ 
        (*3) & f_\tau \textrm{ is minimal.}
      \end{array} \right.
    \end{equation}
  Further, if $f_{\tau_0}$ satisfies
  \begin{equation}
    \label{e.symmetry-addendum}
    f_{\tau_0,\theta+\halb}(-x) \ = \ -f_{\tau_0,\theta}(x) \qquad \forall\thx\in\torus
  \end{equation}
  then $\tau_0\in\Lambda_\infty$.
  \end{thm}
  The proof is given in Section~\ref{ParameterExclusion}.  \smallskip

  \noindent {\em VI. Modified assumptions.} As mentioned in the introduction, it
  is not possible to apply this result to the qpf Arnold circle map due to the
  bounded slope of the fibre maps. In order to make up for this lack of
  hyperbolicity, a particular geometry and symmetry of the forcing has to be
  used. For the twist parameter exclusion carried out here, this is slightly
  more subtle than for the frequency exclusion in \cite{jaeger:2009a} and
  stronger conditions on the forcing are required (see also
  Remark~\ref{r.forcing-structure}).  First, we have to restrict to the case of
  two critical regions with fixed distance $1/2$.
\begin{equation} \label{e.A7'}\textstyle
  \tag{${\cal A}7'$}
 {\cal N} = 2 \ \textrm{ and there holds } \ I^1_0(\tau) = I^2_0(\tau) + \halb  \ .
\end{equation}
Secondly, the slope on the two critical regions must have opposite sign.
\begin{equation}\tag{${\cal A}8'$} \label{e.A8'}
 \partial_\theta f_{\tau,\theta}(x) > s \ \textrm{ on } \ I^1_0\times\kreis \
  \textrm{ and } \partial_\theta f_{\tau,\theta}(x) < -s \ \textrm{ on }
  I^2_0\times\kreis \ .
\end{equation}
Thirdly, as in \cite{jaeger:2009a} we need to ensure that away from the critical
regions the $\theta$-dependence is small. To that end, we suppose ${\cal
I}_0'\ssq \kreis$ is the disjoint union of two open intervals $I^{1'}_0$ and
$I^{2'}_0=I^{1'}_0+\halb$ with $I^k_0\ssq I^{k'}_0$ $(k=1,2)$ and for some $s'
\in (0,s)$ there holds
\begin{equation} \label{eq:refinedbounddth} \tag{${\cal A}9'$}
  |\partial_\theta \fthx| \ < \ s' \ \
  \forall \thx \in (\kreis \smin {\cal I}_0') \times C \ .
\end{equation}
Finally, we need constants $\gamma,L>0$ which provide uniform upper and lower
bounds for the dependence on the twist parameter $\tau$.
\begin{equation} 
   \tag{${\cal A}10'$} \label{e.A10'}
\gamma \ < \ \partial_\tau f_{\tau,\theta}(x) \ < \ L \qquad 
\end{equation}

\begin{thm} \label{t.mr-quantitative2} Let $\omega \in {\cal D}(\sigma,\nu)$ and
  $\delta > 0$. Suppose that $\Lambda \ssq [0,1]$ is an open interval and
  $(F_\tau)_{\Lambda\in[0,1]} \in {\cal P}$ is such that conditions
  (\ref{eq:Diophantine})--(\ref{eq:crossing}) and (\ref{e.A7'})--(\ref{e.A10'})
  are satisfied on $\Lambda$. Let $\eps_0 := \sup_{\iota\in[1,{\cal
      N}],\tau\in\Lambda}|I^\iota_0(\tau)|$. Further, assume there exist
  constants $A,d > 1$ such that
  \begin{eqnarray}
    S & < & A\cdot d \ , \label{e.S<Ad}\\ s & > & d/A \ , \label{e.s>d/A} \\
    \eps_0 & < &A/\sqrt{d} \ , \label{e.eps<1/Ad}\\  \label{e.s'<A} s' & < & A \ .
  \end{eqnarray}

  Then there exist a constant $d_*=d_*(\sigma,\nu,{\cal
    N},\alpha,p,L,\gamma,A,\delta)>0$ with the following property.\smallskip

  If $d > d_*$, then there exists a set
  $\Lambda_\infty=\Lambda_\infty(\omega)\ssq \Lambda$ of measure
  $\Leb(\Lambda_\infty) > \Leb(\Lambda) - \delta$ such that for all $\tau \in
  \Lambda_\infty$ the qpf circle diffeomorphism
  \begin{equation}
    f_\tau \ : \ \thx \mapsto (\thom,\flthx)
  \end{equation}
  satisfies \ref{e.star}. Further, if $f_{\tau_0}$ satisfies
  (\ref{e.symmetry-addendum}) then $\tau_0\in\Lambda_\infty$.
\end{thm}
The proof is given in Section~\ref{RefinedParameterExclusion}. We note that the
precise form of the $d$-dependence in (\ref{e.S<Ad})--(\ref{e.s'<A}) is to some
extent arbitrary and could be stated in a more general way, but we refrain from
introducing even more parameters and refer to Section~\ref{RefinedProof} for
details. The estimates required to deduce
Corollaries~\ref{c.arnold}--\ref{c.tongue-collapse} from this statement will be
carried out in Section~\ref{ArnoldCorollary}.

\section{The basic version of the twist parameter
  exclusion} \label{ParameterExclusion}

\subsection{Critical sets and critical regions.} \label{Previous}

In this section, we will briefly recall the construction from
\cite{jaeger:2009a} and collect the key statements needed for the proof of
Theorem~\ref{t.firstversion}.  The parameter $\tau$, and consequently the map
$f_\tau$, will be fixed.  Nevertheless we keep the dependence on $\tau$ explicit
for the sake of consistency with the later sections. The description of the
dynamics of a suitable qpf circle diffeomorphism $f$ in \cite{jaeger:2009a} is
based on the analysis of certain critical sets ${\cal C}_0 ,{\cal C}_1, {\cal
  C}_2, \ldots$ and critical regions ${\cal I}_0 , {\cal I}_1, {\cal I}_2 ,
\ldots$, which are given as follows.

Given a union ${\cal I}_0(\tau)$ of ${\cal N}$ disjoint open intervals
$I^1_0(\tau) \ldots I^{\cal N}_0(\tau)$ and a monotonically increasing sequence
of integers $\nofolge{M_n}$ with $M_0 \geq 2$, we recursively define
\begin{eqnarray*}
  {\cal A}_n & := & \{\thx \mid \theta \in {\cal I}_n(\tau)-(M_n-1)\omega,\ x
  \in C\} \ , \\ {\cal B}_n & := & \{\thx \mid \theta \in {\cal
  I}_n(\tau)+(M_n+1)\omega,\ x \in E\} \ , \\ {\cal C}_n &:= &
  f_\tau^{M_n-1}({\cal A}_n) \cap f^{-M_n-1}({\cal B}_n) \quad \textrm{and } \\
  \label{eq:defIn}
{\cal I}_{n+1}(\tau) & := & \mathrm{int}(\pi_1({\cal C}_n)) \ .
\end{eqnarray*} 
The crucial observation is the fact that certain {\em `slow
  recurrence' assumptions} on the critical regions ${\cal I}_n$ are
already sufficient to guarantee the nonuniform hyperbolicity of
$f_\tau$. In order to state them, suppose $\nofolge{K_n}$ is a
monotonically increasing sequence of positive integers and
$\nofolge{\eps_n}$ is a non-increasing sequence of positive real
numbers which satisfy $\eps_0 \leq 1$ and $\eps_n \geq 9\eps_{n+1} \
\forall n\in\N_0$. Let
\begin{equation} \nonumber
  {\cal X}_n \ := \ \bigcup_{k=1}^{2K_n M_n} ({\cal I}_n+k\omega) \ \
  \textrm{ and } \ \ {\cal Y}_n \ := \ \bigcup_{j=0}^n \bigcup_{k = -M_j+1}^{M_j+1}
  ({\cal I}_j + k\omega) \ .
\end{equation}
Then the required assumptions on the critical regions are the following.
\begin{equation} \tag*{$({\cal X})_n$} \label{e.Xn} d({\cal I}_j,{\cal X}_j) \ >
  \ 3\eps_j \ \ \ \ \ \forall j = 0 \ld n \ ,
\end{equation}

\vspace{-5ex}
\begin{equation}
  \tag*{$({\cal Y})_n$} \label{e.Yn} d(({\cal I}_j -(M_j-1)\omega) \cup ({\cal
    I}_j+(M_j+1)\omega),{\cal Y}_{j-1}) \ > \ 0 \ \ \ \ \ \forall j=1\ld n \ .
\end{equation}
Let $\beta_0=1$, $\beta_n = \prod_{j=0}^{n-1} \left( 1-\frac{1}{K_j} \right)$
and $\beta=\nLim \beta_n$. Further, define
\begin{equation} \label{e.alpha-infty}
 \alpha_\infty \  = \  \alpha^{2\beta/p - (1-\beta)p} \ .
\end{equation}
\begin{prop}[Propositions 3.10 in \cite{jaeger:2009a}]
  \label{p.sna-existence} If (\ref{eq:Cinvariance})--(\ref{eq:s}) hold,
  $\alpha_\infty > \alpha_1$ and for all $n\in\N_0$ conditions \ref{e.Xn} and
  \ref{e.Yn} are satisfied and ${\cal I}_{n+1}(\tau)=\pi_1({\cal C}_n) \neq
  \emptyset$, then $f_\tau$ has a sink-source orbit and consequently an SNA and an SNR.
\end{prop}
We will also use the following slightly stronger versions of the
above conditions.
\begin{equation} \tag*{$({\cal X'})_n$} \label{e.X'n} d({\cal I}_j,{\cal X}_j) \ >
  \ 9\eps_j \ \ \ \ \ \forall j = 0 \ld n \ ,
\end{equation}

\vspace{-5ex}
\begin{equation}
  \tag*{$({\cal Y'})_n$} \label{e.Y'n} d(({\cal I}_j -(M_j-1)\omega) \cup ({\cal
    I}_j+(M_j+1)\omega),{\cal Y}_{j-1}) \ > \ 2\eps_j \ \ \ \ \ \forall j=1\ld n \ .
\end{equation}

\begin{rem}
  In the proof of Theorems~\ref{t.firstversion} and \ref{t.mr-quantitative2}, we
  will show that for all $\tau\in\Lambda_\infty$ conditions \ref{e.X'n} and
  \ref{e.Y'n} hold. In fact, for our purposes here the weaker conditions
  \ref{e.Xn} and \ref{e.Yn} would be sufficient, since we only need them for the
  application of Proposition~\ref{p.sna-existence}. The reason for using the
  stronger versions \ref{e.X'n} and \ref{e.Y'n} is that we believe these to be
  crucial in a current approach to the proof of Conjecture~\ref{conjecture}, and
  that at the same time this inflicts no extra costs whatsoever.
\end{rem}
Let us briefly review the construction in \cite{jaeger:2009a} that leads to the
statement of Proposition~\ref{p.sna-existence} and provide some more details
that will be used later. Given any point $(\theta_0,x_0) \in \torus$, we denote
its iterates by $(\theta_k,x_k) = f_\tau^k(\theta_0,x_0),\ k\in\Z$.  Now,
(\ref{eq:Cinvariance}) implies that whenever $\theta_0\notin {\cal I}_0(\tau)$
and $x_0\in C$, the forward orbit remains `trapped' in the contracting region
$\kreis \times C$ until $\theta_k$ enters ${\cal I}_0(\tau)$ for the first
time. However, even if $\theta_k \in {\cal I}_0(\tau)$ and the orbit enters the
expanding region at time $k$, that is $x_{k+1} \in E$, it will leave $\kreis
\times E$ again after $M_0$ further iterates unless $\theta_k$ is also contained
in the smaller set ${\cal I}_1(\tau)$. (This is a straightforward consequence of
the definition of ${\cal C}_0$ and ${\cal I}_1(\tau)$.) Following this idea it
is possible, via a purely combinatorial inductive construction, to control the
behaviour of an orbit $(\theta_k,x_k)$ starting in ${\cal A}_{n}$ up to the
first time $\hat k\in\N$ at which $\theta_k$ enters the $(n+1)$-th critical
region ${\cal I}_{n+1}(\tau)$, provided that the `slow recurrence' assumptions
\ref{e.Xn} and \ref{e.Yn} hold \cite[Lemma 3.4]{jaeger:2009a}. In this case, the
finite trajectory will remain in $\kreis \times C$ most of the time and the
fraction of time spent outside this set is at most $1-\beta_n$ \cite[Lemma
3.8]{jaeger:2009a}.

As a consequence, since $M_n\leq \hat k$ by \ref{e.Xn}, it is possible to
control the forward orbit of these points up to time $M_n$, or equivalently the
backward orbit of points in $f_\tau^{M_n}({\cal A}_{n})\supseteq{\cal C}_n$,
which remain trapped in the contracting region most of the time. Similar
findings hold for the forwards orbits of points in $f^{-M_n}_\tau({\cal B}_n)
\supseteq{\cal C}_n$, which remain in the expanding region most of the
time. Combining the combinatorial information about the behaviour of the orbits
with the estimates on the derivatives provided by
(\ref{eq:bounds1})--(\ref{eq:bounds3}) one obtains the following statement.
\begin{lem}[Corollaries 3.7 and 3.9 in \cite{jaeger:2009a}] \label{p.lyaps}
  Suppose (\ref{eq:Cinvariance})--(\ref{eq:bounds3}), \ref{e.Xn} and \ref{e.Yn}
  hold and $\thx \in \closure(f_\tau({\cal C}_n))$. Then for all
  $k\in[0,M_n]$ we have
  \begin{equation}
    \label{e.finite-lyaps}
    \partial_x \flth^{-k}(x) \geq \alpha_\infty^k \quad \textrm{and} 
  \quad \partial \flth^k(x) \geq \alpha_\infty^k \ .
  \end{equation}
  Furthermore, there holds ${\cal C}_0 \supseteq {\cal C}_1 \supseteq \ldots
  \supseteq {\cal C}_{n+1}$ and
\begin{equation}
  \label{e.lastiterates}
  \thx \in ({\cal I}_n(\tau)+\omega) \times E \quad , \quad f_\tau^{-1}\thx \in 
  {\cal I}_n(\tau) \times C \ .
\end{equation}
\end{lem}
This statement rather easily entails Proposition~\ref{p.sna-existence}
(respectively \cite[Proposition 3.10]{jaeger:2009a}): When \ref{e.Xn}
and \ref{e.Yn} hold for all $n\in\N$ and $\thx \in \ncap
\closure\left(f_\tau({\cal C}_n)\right)$, then it follows directly from
(\ref{e.finite-lyaps}) that the point $\thx$ has a positive vertical
Lyapunov exponent both for $f_\tau$ and for its inverse $f^{-1}_\tau$.
Hence, there exist a sink-source-orbit and thus an SNA by
Proposition~\ref{prop:sinksourcesna}.  

Due to Proposition~\ref{p.sna-existence}, the validity of the slow
recurrence conditions \ref{e.Xn} and \ref{e.Yn} together with ${\cal
  I}_n(\tau)\neq\emptyset$ provides a rigorous criterion for the
existence of SNA. The remaining task is to find a positive measure set
of parameters (frequencies $\omega$ in \cite{jaeger:2009a}, parameters
$\tau$ in our setting) for which these assumptions are satisfied for
all $n \in \N$. At this point, a detailed analysis of the geometry of
the critical sets, or more precisely of the sets $f^{M_n}_\tau({\cal
  A}_n)$ and $f^{-M_n}_\tau({\cal B}_n)$ whose intersection equals
$f_\tau({\cal C}_n)$, comes into play. The outcome is that the size of
the connected components of the critical regions ${\cal I}_n(\tau)$
decays super-exponentially with $n$ and that these components depend
`nicely' on the parameter. This will be made precise in
Proposition~\ref{p.lyaps} and Lemma~\ref{l.lambda-dependence} below. 

\begin{figure}[h!]  
\noindent 
\begin{minipage}[t]{\linewidth}  \begin{center}
  \epsfig{file=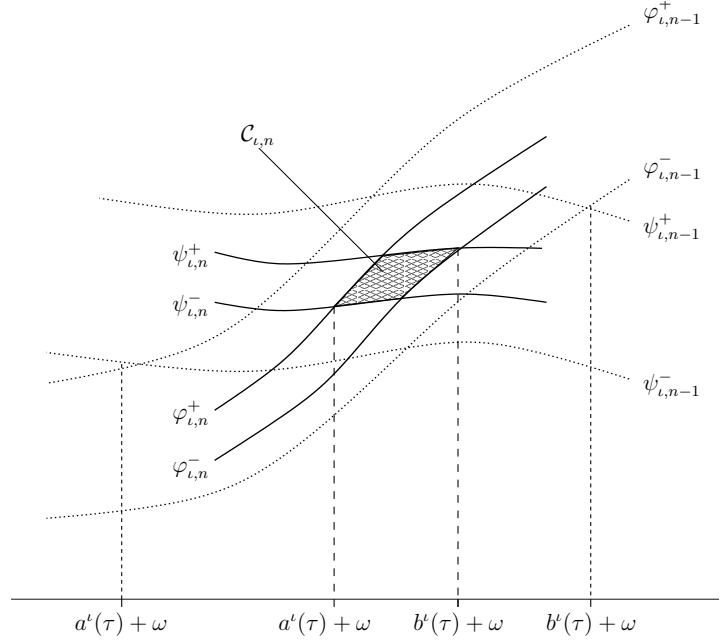, clip=, width=0.7\linewidth}
  \caption{\small The two `strips' $f^{M_n}(A^\iota_n)$ and
    $f^{-M_n}(B^\iota_n)$ intersect each other in a transversal way, producing a
    connected component of ${\cal C}_n$.}\end{center}
         \label{f.crossing}
\end{minipage}
\end{figure}

The main idea behind this geometric analysis is again the fact that for any
connected component $I^\iota_n(\tau)$ of ${\cal I}_n(\tau)$ the first $M_n-1$
forward iterates of the set $A^\iota_n = (I^\iota_n(\tau)-(M_n-1)\omega)\times
C$ remain in the contracting region most of the time. Consequently
$f_\tau^{M_n-1}(A^\iota_n)$ will be a very thin strip, which will moreover be
almost horizontal since the strong contraction `kills' any dependence on
$\theta$. Due to (\ref{eq:s}) the image $f_\tau^{M_n}(A^\iota_n)$ will then have
slope $\geq s/2$. It therefore intersects the image $f_\tau^{-M_n}(B^\iota_n)$
of $B^\iota_n = (I^\iota_n(\tau)+(M_n+1)\omega)\times E$ in a transversal way,
since this set is a very thin horizontal strip by the same reasoning as for
$f_\tau^{M_n-1}(A^\iota_n)$. (Note that the expanding region $\kreis\times E$ is
contracted by the inverse $f^{-1}$.) Consequently the resulting intersection
$f_\tau^{M_n}(A^\iota_n) \cap f_\tau^{-M_n}(B^\iota_n)$ always has the geometry
depicted in Figure~\ref{f.crossing} and projects to a very small interval
$I^\iota_{n+1}(\tau)+\omega$.

In order to give a more detailed quantitative version of this heuristic
description, we first define the `bounding graphs' of the sets
$f_\tau^{M_n}(A^\iota_n)$ and $f_\tau^{-M_n}(B^\iota_n)$.  For $\theta \in
I^\iota_n(\tau)+\omega$, let
\begin{equation} \label{eq:phipsidef}
  \varphi_{\iota,n}^\pm(\theta,\tau) \ := \
  f_{\tau,\theta-M_n\omega}^{M_n}(c^\pm) \quad \textrm{ and } \quad
  \psi_{\iota,n}^\pm(\theta,\tau) \ := \
  f_{\tau,\theta+M_n\omega}^{-M_n}(e^\pm)\ .
\end{equation}
Note that we have
\begin{eqnarray} \label{e.phin-characterisation} f_\tau^{M_n}(A^\iota_n) & =
  & \{ \thx \mid \theta \in I_n^\iota(\tau)+\omega,\ x \in
  [\varphi_{\iota,n}^-(\theta,\tau) ,\varphi_{\iota,n}^+(\theta,\tau)] \}
  \ , \\ \label{e.psin-characterisation} f_\tau^{-M_n}(B^\iota_n) & = & \{
  \thx \mid \theta \in I_n(\tau)+\omega,\ x \in
  [\psi_{\iota,n}^-(\theta,\tau), \psi_{\iota,n}^+(\theta,\tau)] \} \ .
  \end{eqnarray}

Using this notation, we can now restate the following estimates from
\cite{jaeger:2009a}.

\begin{prop}[Proposition~3.11 and Lemma~3.14 in \cite{jaeger:2009a}]
  \label{p.In-size} Suppose
  (\ref{eq:Cinvariance})--(\ref{eq:s}) hold and \ref{e.Xn} and
  \ref{e.Yn} are satisfied.  Then there exists a constant $\alpha_1 =
  \alpha_1(s,S)$ such that the following holds: If $\alpha_\infty >
  \alpha_1$ then \romanlist
\item ${\cal I}_n(\tau)$ and ${\cal I}_{n+1}(\tau)$ consist of exactly
  ${\cal N}$ connected components and each component of ${\cal I}_n(\tau)$
  contains exactly one component of ${\cal I}_{n+1}(\tau)$.
\item For all connected components of $I^\iota_{n+1}(\tau)$ of ${\cal
    I}_{n+1}(\tau)$ there holds
  \begin{equation} \label{e.In-size}
    |I^\iota_{n+1}(\tau)| \ \leq \ 2\alpha_\infty^{-M_n}/s \ .
  \end{equation}
\item Either 
  \begin{eqnarray}\label{e.transversal-intersection}
\inf_{\theta\in{\cal I}_n(\tau)+\omega}
\partial_\theta\varphi^\xi_{\iota,n}(\theta,\tau) - \partial_\theta\psi^\zeta_{\iota,n}(\theta,\tau) \
  > & 0 & \quad \forall \xi,\zeta \in \{+,-\} \qquad \textrm{or} \\\label{e.downwards}
\sup_{\theta\in{\cal I}_n(\tau)+\omega}
\partial_\theta\varphi^\xi_{\iota,n}(\theta,\tau) - \partial_\theta\psi^\zeta_{\iota,n}(\theta,\tau) &
  < & 0 \quad \forall \xi,\zeta \in \{+,-\} \ .
\end{eqnarray}
\listend
\end{prop}
We note that part (iii) can be seen justification for the picture in
Figure~\ref{f.crossing}, insofar as
(\ref{e.transversal-intersection}), respectively (\ref{e.downwards})
ensures that each of the pairs of curves $\varphi^\xi_{\iota,n}$ and
$\psi^\zeta_{\iota,n}$ intersect in exactly one point and
$f_\tau^{M_n}(A^\iota_n)$ crosses $f^{-M_n}(B^\iota_n)$ either upwards
(as depicted) or downwards.\smallskip

It remains to obtain a good control on the dependence of the critical regions
${\cal I}_n(\tau)$ on $\tau$. This will be the content of the next section. On
the technical level, this is the crucial difference in comparison to the
construction in \cite{jaeger:2009a}.

\subsection{Dependence of the critical regions on $\tau$.} \label{LambdaDependence}

In order to perform the parameter exclusion with respect to $\tau$, we
need to show that different connected components of ${\cal I}_n(\tau)$
move with different speed as the parameter $\tau$ changes. This is
ensured by the following lemma.
\begin{lem}
  \label{l.lambda-dependence} Suppose
  (\ref{eq:Cinvariance})--(\ref{e.d}) hold and \ref{e.Xn} and
  \ref{e.Yn} are satisfied.  Then there exists a constant $\alpha_2 =
  \alpha_2(s,S,L,\eta)$ such that the following holds:

  If $\alpha_\infty > \alpha_2$ then all connected components of ${\cal
    I}_{n+1}(\tau)$ are differentiable in $\tau$. Further,
  \begin{eqnarray}
    \label{e.lambda-dependence}
    |\partial_\tau I^\iota_{n+1}(\tau)| & < & 2L/s \qquad 
    \forall \iota \in \{1\ld{\cal N}\} \qquad \textrm{ and } \\\label{e.dlambda-bound}
    |D_\tau(I^\iota_{n+1}(\tau),I^\kappa_{n+1}(\tau))| & > & \frac{\eta}{2} \qquad \qquad
   \forall \iota\neq\kappa\in\{1\ld {\cal N}\} \ . 
  \end{eqnarray}
\end{lem}
\myproof We write
$I_{n+1}^\iota(\tau)=(a^\iota(\tau),b^\iota(\tau)$ and
show that
\begin{equation}
  \label{e.endpoint-ldep}
  \partial_\tau a^\iota(\tau),\ \partial_\tau b^\iota(\tau)\ \in  
  \ B_\Delta(Q^\iota(\tau)) \quad \forall \iota\in[1,{\cal N}]) 
\end{equation}
where $\Delta:=\min\{\eta/4,L/s\}$.  Due to (\ref{e.d}) this immediately implies
(\ref{e.dlambda-bound}). Further, since $Q^\iota(\tau)\ssq [-L/s,L/s]$ by
(\ref{eq:s}) and (\ref{eq:bounddlambda}) we also obtain~(\ref{e.lambda-dependence}).

We carry out the proof only for $\partial_\tau a^\iota(\tau)$ since
the other endpoint can be treated in the same way. Similarly, we
assume that the crossing between $f^{M_n}_\tau(A^\iota_n)$ and
$f^{-M_n}_\tau(B^\iota_n)$ is upwards, that is, case
(\ref{e.transversal-intersection}) in Proposition~\ref{p.In-size}(iii)
holds. Again, the other case can be treated similarly.  Then, as can be
seen from the picture in Figure~\ref{f.crossing}, $a^\iota(\tau)$ is
characterized by the equality
\begin{equation} \nonumber
  \varphi^+_{\iota,n}(a^\iota(\tau)+\omega,\tau) - 
  \psi^-_{\iota,n}(a^\iota(\tau)+\omega,\tau) \ = \ 0 \ . 
\end{equation}
Application of the Implicit Function Theorem therefore yields 
\begin{equation} \label{e.IFT-Formula}
  \partial_\tau a^\iota(\tau) \ = \ 
  - \ \frac{\partial_\tau\left(\varphi^+_{\iota,n}- 
      \psi^-_{\iota,n}\right)(a(\tau)+\omega,\tau)}{\partial_\theta \left(\varphi^+_{\iota,n}- 
      \psi^-_{\iota,n}\right)(a(\tau)+\omega,\tau)} \ .
\end{equation}
We start by deriving an estimate on the numerator. Let $\theta :=
a^\iota(\tau)+\omega$ and $x:= \varphi^+_{\iota,n}(\theta,\tau) =
f^{M_n}_{\tau,\theta-M_n\omega}(c^+)$.  Let $\theta_0 =
\theta-M_n\omega$ and $\theta_k = \theta_0+k\omega$. Further, let
$x_0=c^+$ and $x_k=f^k_{\tau,\theta_0}(x_0)$. Note that thus $\thx =
(\theta_{M_n},x_{M_n})$. Differentiating with respect to $\tau$ we obtain
\begin{equation}
  \partial_\tau \varphi^+_{\iota,n}(\theta,\tau) \ = \ 
  \partial_\tau f_{\tau,\theta_{M_n-1}} 
  (x_{M_n-1})  \label{e.rqiota}
  \  + \  \ \underbrace{\sum_{k=1}^{M_n-1} 
    \partial_x f^{M_n-k}_{\tau,\theta_k}(x_k)
    \cdot \partial_\tau f_{\tau,\theta_{k-1}} 
    (x_{k-1})}_{=:r_1^\iota} \ .
\end{equation}
Since $
\partial_x f^{M_n-k}_{\tau,\theta_{k}}(x_k) = \left( \partial_x
  f_{\tau,\theta}^{-M_n+k}(x)\right)^{-1} $ and $(\theta,x) \in
\closure\left(f_\tau({\cal C}_n)\right)$, the estimate
(\ref{e.finite-lyaps}) in Lemma~\ref{p.lyaps} yields $\partial_x
f^{M_n-k}_{\tau,\theta_k}(x_k) \leq \alpha_\infty^{-k}$.  Together
with the upper bound $L$ on the derivative with respect to $\tau$
provided by (\ref{eq:bounddlambda}) we obtain $|r_1^\iota| \ \leq \ L
\cdot \sum_{k=1}^\infty \alpha_\infty^{-k}$. 

Now let $\theta=a^\iota(\tau)+\omega$ as before,
$\xi=\psi^-_{\iota,n}(\theta,\tau)$, $\vartheta_0=\theta+M_n\omega$,
$\vartheta_k=\vartheta_0+k\omega$,$\xi_0=e^-$ and
$\xi_k=f^k_{\tau,\vartheta_0}(\xi_0)$. Note that $(\vartheta_{-M_n},xi_{-M_n}) =
(\theta,\xi)$.  Differentiating with respect to $\tau$ yields
\begin{eqnarray} \nonumber r^\iota_2 & := & \partial_\tau
  \psi^-_{\iota,n}(\theta,\tau) \  = \ 
    \sum_{k=1}^{M_n} \partial_x f_{\tau,\vartheta_{-k}}^{-M_n+k} 
  (\xi_{-k}) \cdot 
  \partial_\tau f^{-1}_{\tau,\vartheta_{-k+1}}
  (\xi_{-k+1})\\ \label{e.psicontrol}
  & = & -\ \sum_{k=1}^{M_n} \left(\partial_x f^{M_n-k}_{\tau,\vartheta_{-M_n}}
    (\xi_{-M_n}) \right)^{-1}\cdot 
  \left(\partial_x f_{\tau,\vartheta_{-k}}(\xi_{-k})\right)^{-1} 
  \cdot \partial_\tau f_{\tau,\vartheta_{-k}}(\xi_{-k}) \\
  & = & -\ \sum_{k=1}^{M_n} \left(\partial_x 
  f^{M_n-k+1}_{\tau,\vartheta_{-M_n}}(\xi_{-M_n})\right)^{-1} \cdot 
  \partial_\tau f_{\tau,\vartheta_{-k}}(\xi_{-k}) \ . \nonumber
\end{eqnarray}
Using (\ref{eq:bounddlambda}) and (\ref{e.finite-lyaps}) again we
obtain $|r^\iota_2| \leq L \cdot \sum_{k=1}^\infty
\alpha_\infty^{-k}$.  If we let $r^\iota = r^\iota_1 - r^\iota_2$ and
use that $\theta_{M_n-1}=\theta-\omega$ and $x_{M_n-1} =
f_{\tau,\theta}^{-1}(x)$ then
\begin{eqnarray}
  \partial_\tau\left(\varphi^+_{\iota,n}- 
    \psi^-_{\iota,n}\right)(\theta,\tau) & = & \partial_\tau 
  f_{\tau,\theta-\omega}\left(f_{\tau,\theta}^{-1}(x)\right) \ + \ r^\iota \\
  \textrm{with } \quad |r^\iota|  & \leq & 2L/(\alpha_\infty-1) \ .
\end{eqnarray}
If we replace $\partial_\tau$ by $\partial_\theta$ in these computations and
use (\ref{eq:bounddth}) instead of (\ref{eq:bounddlambda}), then we obtain in
exactly the same way that
\begin{eqnarray}
  \partial_\theta\left(\varphi^+_{\iota,n}- 
    \psi^-_{\iota,n}\right)(\theta,\tau) & = & \partial_\theta 
  f_{\tau,\theta-\omega}\left(f_{\tau,\theta}^{-1}(x)\right) \ + \ q^\iota \\
\textrm{with } \quad |q^\iota| &  \leq &  2S/(\alpha_\infty-1) \ .
\end{eqnarray}
Now $\theta-\omega \in {\cal I}_n \ssq {\cal I}_0$, such that
$\left|\partial_\theta
  f_{\tau,\theta-\omega}\left(f_{\tau,\theta}^{-1}(x)\right)\right| > s$
by (\ref{eq:s}) and
\[
\partial_\tau
  f_{\tau,\theta-\omega}\left(f_{\tau,\theta}^{-1}(x)\right)/\partial_\theta
  f_{\tau,\theta-\omega}\left(f_{\tau,\theta}^{-1}(x)\right) \ \in \
Q^\iota(\tau)
\]
by the definition of $Q^\iota(\tau)$. Furthermore, it follows from the above
estimates that $|r^\iota|$ and $|q^\iota|$ go to zero as $\alpha_\infty \to
\infty$. Hence, for sufficiently large $\alpha_\infty$ we have
\begin{equation} \label{e.lamda-dep-form}
  \partial_\tau a^\iota(\tau) \ = 
  \ - \ \frac{\partial_\tau f_{\tau,\theta-\omega}\left(f_{\tau,\theta}^{-1}(x)\right) + r^\iota}{
    \partial_\theta f_{\tau,\theta-\omega}\left(f_{\tau,\theta}^{-1}(x)\right) + q^\iota} \
  \in \ B_{\Delta}(Q^\iota(\tau)) \ .
\end{equation}
Furthermore, it can be seen from the above estimates that the
largeness condition on $\alpha_\infty$ only depends on the constants
$s,S,L$ and $\eta$.  \qed

\subsection{Preliminaries for the parameter
  exclusion} \label{ParameterExclusionPreliminaries} We now collect some
preliminary statements for the parameter exclusion. The setting is an abstract
one that does not depend on the previous dynamical construction. We first fix an
integer ${\cal N}$ and sequences $\nofolge{K_n}$ and $\nofolge{\eps_n}$ with the
same properties as in Section~\ref{Previous} and a sequence $\nofolge{N_n}$ of
positive integers that satisfy
\begin{equation}
  \label{e.N}\tag{${\cal N}1$} N_0 \geq 3 \quad \textrm{and} \quad N_{n+1} 
  \geq 2K_nN_n \ \forall n\in\N_0 \ .
\end{equation}
We denote by ${\cal S}(\kreis)$ the set of all subsets of $\kreis$. Let
$\Lambda\ssq[0,1]$ be an open interval. Then we simply assume that we are given
a sequence of mappings
\[ {\cal I}_n \ : \quad \Lambda \times \N^n \to {\cal S}(\kreis) \ , \quad
(\tau,M_0\ld M_{n-1}) \ \mapsto \  {\cal I}_n(\tau) = {\cal I}_n(\tau,M_0\ld
M_{n-1}) \ 
\]
The dependence of ${\cal I}_n(\tau)$ on $M_0 \ld M_{n-1}$ will be kept
implicit. We let
\[
{\cal P}_n = {\cal P}_n(M_0\ld M_n) \ := \ \{ \tau \in \Lambda \mid
 ({\cal X'})_n \textrm{ and } ({\cal Y'})_n  \textrm{ hold} \} \ .
\]
Here $({\cal X'})_n$ and $({\cal Y'})_n$ are understood
as conditions on the sets ${\cal I}_j(\tau)$ ($j\in[0,n]$) for fixed
$\omega\in\kreis$. Furthermore, we assume that the following
conditions are satisfied.
\begin{equation}
  \label{e.P} \tag{$\cal P$}
\left\{ \quad 
  \begin{array}{l} \textrm{Suppose }  M_0\ld M_n \textrm{ with } M_i \in [N_i,2K_iN_i)
    \ \forall i=[0,n] \textrm{ are fixed.} \\  
    \textrm{Then for all } \tau \in {\cal P}_n(M_0\ld M_n) \textrm{ and } 
    j\in[0,n+1] \textrm{ there holds } \\ \ \\
  \begin{array}{cl}
    ({\cal P}1) & {\cal I}_j(\tau) \textrm{ is open and consists of exactly } {\cal N}
    \textrm{ connected components }\\ & I^1_j(\tau) \ld I^{\cal N}_j(\tau).
    \textrm{ If } j \leq n \textrm{ then } I^\iota_{j+1}(\tau) \ssq I^\iota_j(\tau) \ 
    \forall \iota \in [1,{\cal N}]. \\ \ \\
    ({\cal P}2) & \textrm{The set } {\cal P}_j(M_0\ld M_j) 
    \textrm{ is open and all } I^\iota_j\textrm{ with }
    j\in [1,{\cal N}] \\ &
    \textrm{ are differentiable w.r.t.\ } \tau \textrm{ on } {\cal P}_j(M_0\ld M_j);
    \\ \ \\   ({\cal P}3) & |I^\iota_j(\tau)| \leq \eps_j \quad \forall 
    \iota \in [1,{\cal N}] \ ; \\ \ \\
    ({\cal P}4) & |\partial_\tau I^\iota_j(\tau)| \leq 2L/s \quad 
    \forall \iota \in [1,{\cal N}] \ ;\\ \ \\
    ({\cal P}5) & |D_\tau(I^\iota_j(\tau),I^\kappa_j(\tau))| \geq \eta/2 \quad 
    \forall \iota \neq \kappa \in [1,{\cal N}] \ ;
  \end{array}
\end{array}
\right.
\end{equation}
\begin{rem} \label{b.parameter-exclusion-assumption} Note that if
  (\ref{eq:Cinvariance})--(\ref{e.d}) are satisfied and $\alpha_\infty$ is
  sufficiently large, then the fact that (\ref{e.P}) holds for the sets ${\cal
    I}_j(\tau)$ defined by (\ref{eq:defIn}) is exactly the content of the
  previous sections. $({\cal P}1)$ follows by induction from
  Proposition~\ref{p.In-size}.  Lemma~\ref{l.lambda-dependence} implies that the
  connected components of ${\cal I}_j$ are differentiable with respect to
  $\tau$, which in turn yields the openness of the conditions $({\cal X}')_j$ and
  $({\cal Y'})_j$ and hence of ${\cal P}_j(M_0\ld M_j)$, such that $({\cal P}2)$
  holds as well. $({\cal P}3)$ follows from Proposition~\ref{p.In-size} and
  finally $({\cal P}4)$ and $({\cal P}5)$ are again a consequence of
  Lemma~\ref{l.lambda-dependence}.
\end{rem}

In each step of the parameter exclusion we will have to ensure that the set
${\cal P}_n \smin {\cal P}_{n+1}$ of excluded parameters is small. In other
words, we have to show that for most $\tau \in {\cal P}_n$ the conditions
$({\cal X}')_{n+1}$ and $({\cal Y}')_{n+1}$ are satisfied for a suitable
$M_{n+1}$ (that we allow to depend on $\tau$). This is greatly
simplified by the fact that $({\cal Y})_{n+1}$ `comes for free'.
\begin{lem}[Lemma 3.16 in \cite{jaeger:2009a}]
  \label{l.Yn+1-Mexists} Suppose that $M_0 \ld M_n$ with $M_j \in [N_j,2N_j) \
  \forall j\in[0,n]$ are fixed. Further, assume that (\ref{e.N}), $({\cal P}1)$
  and $({\cal P}3)$ hold and
  \begin{equation}
    \label{e.K}\tag{${\cal K}$}
    \sum_{j=0}^\infty \frac{1}{K_j} \ < \ \frac{1}{6{\cal N}^2} \ .
  \end{equation}
  Then for all $\tau \in {\cal P}_n(M_0\ld M_n)$ there exists an integer $M(\tau)
  \in [N_{n+1},2N_{n+1})$ such that
  \begin{equation}
    \label{e.strong-Yn+1}
    d\left(({\cal I}_{n+1}(\tau)-(M(\tau)-1)\omega)\cup({\cal I}_{n+1}+(M(\tau)+1)\omega)\ , 
      \ {\cal Y}_n\right) \ > \ 3\eps_n \ . 
  \end{equation}
\end{lem}
We remark that the version of this lemma in \cite{jaeger:2009a} actually
contains some additional assumptions, but these are not used in the proof. (For
the sake of brevity, the standard hypothesis were assumed throughout the
respective section in \cite{jaeger:2009a}.) In order to obtain an estimate on
the set of $\tau\in {\cal P}_n$ that do not satisfy $({\cal X})_{n+1}$, the
following lemma is needed.
\begin{lem}\label{l.pex-basic}
  Suppose $\Lambda \ssq [0,1]$ is an interval and ${\cal I} : \Lambda \to {\cal
    S}(\kreis)$ is such that for all $\tau\in\Lambda$ the set ${\cal
    I}(\tau) \ssq \kreis$ consists of ${\cal N}$ connected components
  $I^1(\tau) \ld I^{\cal N}(\tau)$ of length $|I^\iota(\tau)| \leq
  \delta$ which satisfy
  \begin{equation} \label{e.pex-basic1}
  |D_\tau(I^\iota(\tau),I^\kappa(\tau))| \geq \eta/2 \quad
 \forall \iota \neq \kappa \in [1,{\cal N}] \ . 
  \end{equation}
  Further, assume that 
  \begin{equation} \label{e.pex-diophantine}
    d(I^\iota(\tau),I^\iota(\tau)+n\omega) \ > \ \eps \quad \forall
    \tau \in \Lambda,\ n \in [1,M],\ \iota \in[1,{\cal N}] \ . 
  \end{equation}
  Then the set
  \[
  \Upsilon \ := \ \left\{ \tau \in \Lambda \ \left| \ d\left({\cal
          I}(\tau),\bigcup_{j=1}^M {\cal I}(\tau)+n\omega\right) \leq \eps
    \right.\right\}
  \]
  has measure $\leq 8{\cal N}^2M\frac{\delta+\eps}{\eta}$ and consists of at
  most $2{\cal N}^2M-1$ connected components.
\end{lem}
\myproof Fix $\iota\neq \kappa \in [1,{\cal N}]$ and $n\in[1,M]$. As
$I^{\iota}(\tau)$ and $I^\kappa(\tau)$ are disjoint for all $\tau \in
\Lambda$ and due to (\ref{e.pex-basic1}), the set of $\tau$ with
$d\left(I^\iota(\tau),(I^\kappa(\tau)+n\omega)\right)\leq \eps$ consists of at most
two intervals of length $\leq 4(\delta+\eps)/\eta$. Summing up over all
$\iota,\kappa$ and $n$ yields the statement. \qed
\medskip

For any $n\geq 1$, let
\begin{eqnarray}
  v_n & = & \ 32{\cal N}^2 K_{n+1}N_{n+1}\cdot \frac{L}{s\eps_{n-1}} \quad \textrm{and}\\
  u_n & = & 1280{\cal N}^2 K_{n+1}N_{n+1}\cdot \frac{L }{s\eta}  \cdot \frac{\eps_{n}}{\eps_{n-1}} \ .
\end{eqnarray}
Further, let $v_0 = 8{\cal N}^2K_0N_0$ and $u_0 = 320{\cal
  N}^2K_0N_0\eps_0/\eta$.
\begin{lem} \label{l.pex-component}
  Suppose $\omega \in {\cal D}(\sigma,\nu)$, (\ref{e.N}), (\ref{e.K}) and
  (\ref{e.P}) hold and 
  \begin{equation}
    \label{e.N2} \tag{${\cal N}2$}
    10\eps_n \ < \ \sigma\cdot(4K_nN_n)^{-\nu} \quad \forall n\in\N_0 \ . 
  \end{equation}
  Further, assume that $M_0 \ld M_n$ with $M_j \in [N_j,2N_j) \ \forall
  j\in[0,n]$ are fixed and $\Gamma \ssq {\cal P}_n(M_0\ld M_n)$ is an
  interval. Then for some $r\leq v_{n+1}$ there exist disjoint intervals
  $\Gamma_1 \ld \Gamma_r \ssq \Gamma$ and numbers $M^k \in [N_{n+1},2N_{n+1}),\
  k\in[1,r]$ such that
  \begin{eqnarray}
    \Gamma^k & \ssq & {\cal P}_{n+1}(M_0\ld M_n,M^k) \quad \textrm{and} \\
    \quad \sum_{k=1}^r \Leb(\Gamma^k) & \geq & \Leb(\Gamma) - u_{n+1} \ .
  \end{eqnarray}
\end{lem}
\myproof Divide $\Gamma$ into at most $\frac{4L}{s\eps_n}$ intervals $\Omega_i$
of length of length $\leq \frac{s \eps_n}{2L}$. Denote the midpoint of
$\Omega_i$ by $\tau_i$ and choose $\tilde M^i=M(\tau_i)$ according to
Lemma~\ref{l.Yn+1-Mexists} such that (\ref{e.strong-Yn+1}) holds for
$\tau_i$. Then due to $({\cal P}4)$ we obtain that $({\cal Y})_{n+1}$ holds for
all $\tau \in \Omega^i$. Application of Lemma~\ref{l.pex-basic} with
$M=2K_{n+1}\tilde M^i < 4K_{n+1}N_{n+1}$, $\delta = \eps_{n+1}$ and
$\eps=9\eps_{n+1}$ yields the existence of a set $\tilde \Omega_i \ssq {\cal
  P}_{n+1}(M_0\ld M_n,\tilde M^i)$ of measure $\geq \Leb(\Omega_i)-320{\cal
  N}^2K_{n+1}N_{n+1}\cdot \frac{\eps_{n+1}}{\eta}$ and with at most $8{\cal
  N}^2K_{n+1}N_{n+1}$ connected components. Note that the fact that
(\ref{e.pex-diophantine}) holds follows from the Diophantine condition on
$\omega$ together with (\ref{e.N2}) and $({\cal P}3)$. Relabelling the connected
components of the sets $\tilde\Omega_i$ and summing up over all $i$ yields the
statement.  \qed \medskip

Let $V_{-1}=1$ and $V_n \ = \ \prod_{j=0}^n v_j$ for $n\geq 0$. 
\begin{prop}\label{p.pex}
  Suppose $\omega\in{\cal D}(\sigma,\nu)$, $({\cal N}1$-$2)$, (\ref{e.K}) and
  (\ref{e.P}) hold and
  \begin{equation}
    \label{e.gamma}
    m \ := \ \Leb(\Lambda) - \sum_{n=0}^\infty V_{n-1} u_n . 
  \end{equation}
  \alphlist 
  \item
  Then there exists a set $\Lambda_\infty \ssq \Lambda$ of measure $\geq m$
  with the following property: For all $\tau \in \Lambda_\infty$ there exists
  a sequence $\nofolge{M_n(\tau)}$ with $M_n(\tau)\in[N_n,2N_n) \ \forall
  n\in\N$ such that $\tau \in \bigcap_{n\in\N} {\cal P}_n(M_0(\tau)\ld
  M_n(\tau))$.
\item If there exists $M_0\in[N_0,2N_0)$ with ${\cal
    P}_0(N_0)=\Lambda$, then $\Lambda_\infty$  can be chosen with measure $\geq m+u_0$.
\listend
\end{prop}
\myproof We construct a nested sequence of sets $\Lambda_n$ with the following
properties: 
\romanlist
\item $\Lambda_n$ consists of $\rho_n \leq V_n$ disjoint intervals $\Lambda^1_n
  \ld \Lambda^{\rho_n}_n$;
\item $\Leb(\Lambda_n) \geq \Leb(\Lambda)-\sum_{i=0}^n V_{n-1}u_n$;
\item For each $i\in[1,\rho_n]$ there exist numbers $M_0^{n,i} \ld M_n^{n,i}$
  such that $\Lambda_n^i \ssq {\cal P}_n(M_0^{n,i}\ld M_n^{n,i})$;
\item For each $k\leq n$ and each $i\in[1,\rho_n]$ there exists a unique $\kappa
  \in[1,\rho_k]$ such that $\Lambda^i_n \ssq \Lambda^\kappa_k$ and
  $M^{n,i}_j=M^{k,\kappa}_j \ \forall j\in[0,k]$. \listend 

  The set $\Lambda_\infty = \bigcap_{n\in\N_0} \Lambda_n$ then clearly has the
  properties required in (a), and for (b) it suffices to note that if ${\cal
    P}(M_0)=\Lambda$, then obviously a measure of $u_0$ is gained in the first
  step of the construction.\smallskip

\noindent
For $n=0$ we choose $M_0 \in [N_0,2N_0)$ arbitrarily and let $\Lambda_0={\cal
  P}_0(M_0)$. The fact that it has the required properties follows directly from
Lemma~\ref{l.pex-basic}. Now suppose that $\Lambda_0 \ld \Lambda_n$ with the
above properties exist. Then for each $i\in[1,\rho_n]$ we can apply
Lemma~\ref{l.pex-component} and obtain a union of at most $v_{n+1}$ intervals
with overall measure $\geq \Leb(\Lambda_n^i) - u_{n+1}$. Doing this for the at
most $V_n$ components of $\Lambda_n$ yields the required set
$\Lambda_{n+1}$.
 \qed

\subsection{Minimality and the uniqueness of SNA.}

As a first step in the proof of Theorem~\ref{t.firstversion} below, we will
define the set $\Lambda_\infty$ and show that for all $\tau \in \Lambda_\infty$
the slow-recurrence conditions \ref{e.X'n} and \ref{e.Y'n} hold. Once this is
accomplished, the parameter dependence on $\tau$ does not play a role anymore
and we can consider the map $f_\tau$ as being fixed. The existence of an SNA and
an SNR then follows from Proposition~\ref{p.sna-existence}, and it
remains to prove the uniqueness and one-valuedness of the invariant graphs and
the minimality of $f$. However, this second step has already been carried out in
\cite{jaeger:2009a} and the proof given there literally remains true in our
setting. Instead of repeating it here, we just give a precise formulation of the
formal statement that can be deduced from \cite{jaeger:2009a}.
\begin{prop} \label{p.minimality}
  Suppose $f_\tau$ satisfies (\ref{eq:Cinvariance})--(\ref{eq:s}), \ref{e.Xn} and
  \ref{e.Yn} hold for all $n\in\N$ and 
\begin{equation} \label{e.minimality-condition}
\Leb\left(\bigcup_{n=0}^{\infty} \bigcup_{k=-M_n-1}^{M_n+1} {\cal
    I}_n-k\omega \right) \ < \ \frac{1}{4+4p^2}\ .
\end{equation}
Then $f_\tau$ has a unique SNA and SNR which are both one-valued. Further the
dynamics are minimal.
\end{prop}
\myproof See Sections 3.6 and 3.7 in \cite{jaeger:2009a}. \qed

\subsection{Proof of Theorem~\ref{t.firstversion}.}

Fix an integer $t\geq 4$ such that $2^{-t+2}/{\cal N}^2 \leq
\log((p^2+2)/(p^2+1))$ and let $K_n = 2^{n+t}{\cal N}^2$. Then it is easy to
check that (\ref{e.K}) is verified and furthermore $\beta$ in
(\ref{e.alpha-infty}) is larger than $(p^2+1)/(p^2+2)$. This in turn implies
that $\alpha_\infty$ defined by (\ref{e.alpha-infty}) is larger than
$\alpha^{1/p}$. Suppose that (\ref{eq:Cinvariance})--(\ref{e.d}) hold and
$\alpha^{1/p} > \max\{\alpha_1,\alpha_2\}$, where $\alpha_1$ and $\alpha_2$ are
the constants from Proposition~\ref{p.In-size} and
Lemma~\ref{l.lambda-dependence}. Then as mentioned in
Remark~\ref{b.parameter-exclusion-assumption}, the critical regions ${\cal I}_n$
defined dynamically by (\ref{eq:defIn}) satisfy (\ref{e.P}) when viewed as
mappings ${\cal I}_n : \Lambda \times \N^n \to {\cal S}(\kreis)$. In order to
determine the set $\Lambda_\infty$ by applying Proposition~\ref{p.pex}, it only
remains to show that by an appropriate choice of the sequences
$\nofolge{\eps_n}$ and $\nofolge{N_n}$ we can ensure that (\ref{e.N}) and
(\ref{e.N2}) hold and that the sum $\sum_{n=0}^\infty V_{n-1} u_n$ in
(\ref{e.gamma}) is smaller than $\delta$.

In order to do so, we let $N_0 = 3$ and $N_{n+1} = \alpha^{N_n/qp}$, where
$q=\max\{8,2\nu\}$. Further, we let $\eps_0=\sup_{\iota\in[1,{\cal
    N}],\tau\in\Lambda} |I^\iota_0(\tau)|$ and $\eps_{n+1} = \frac{2}{s} \cdot
\alpha^{-N_n/p}$.  Note that these sequences grow, respectively decay,
super-exponentially. Therefore it is easy to see that with this choice
(\ref{e.N}) and (\ref{e.N2}) are satisfied for sufficiently large $\alpha$ and
sufficiently small $\eps_0$.  In the following estimates we assume that $\alpha$
is chosen sufficiently large and indicate the steps in which this fact is used
by placing $(\alpha)$ over the respective inequality signs. For any $n\in\N_0$
we have
\begin{eqnarray*}
  v_{n+1}& = & 32{\cal N}^2K_{n+1}N_{n+1} \cdot \frac{L}{s\eps_n}  \\
  & = & 16\cdot 2^{n+t}{\cal N}^2  L \cdot \alpha^{N_{n}/qp+N_{n-1}/p} 
  \ \stackrel{(\alpha)}{\leq} \ \alpha^{N_n/4} \quad \textrm{and} \\
  u_{n+1} & = & 1280{\cal N}^2K_{n+1}N_{n+1}\cdot\frac{L}{s\eta}\cdot \frac{\eps_{n+1}}{\eps_n} \\
  & = & 1280\cdot 2^{n+t}{\cal N}^2\cdot\frac{L}{s\eta}\cdot \alpha^{N_n/qp+N_{n-1}/p-N_n/p} 
    \ \stackrel{(\alpha)}{\leq} \ \alpha^{-3N_n/4p} \ .
\end{eqnarray*}
By induction, we obtain that $V_n = \prod_{j=0}^n v_j \stackrel{(\alpha)}{\leq}
\alpha^{N_n/4}$ (note that $v_0 = 8{\cal N}^2K_0N_0 \stackrel{(\alpha)}{\leq}
\alpha^{N_0/4}$). Altogether, this yields that $m$ in Proposition~\ref{p.pex}
satisfies $m \geq \Leb(\Lambda) - u_0 - \sum_{n=0}^{\infty} \alpha^{-N_n/2p}$. 
As $u_0 = 320{\cal N}^2K_0N_0\eps_0/\eta$, this lower bound goes to
$\Leb(\Lambda)$ as $\eps_0 \to 0$ and $\alpha \to \infty$.

Hence, Proposition~\ref{p.pex} yields the existence of a set $\Lambda_\infty\ssq
\Lambda$ of measure $\Leb(\Lambda_\infty) > \Leb(\Lambda) - \delta$ such that
for all $\tau \in \Lambda_\infty$ the conditions \ref{e.X'n} and \ref{e.Y'n} are
satisfied. Fix $\tau\in\Lambda_\infty$.  As $\alpha_\infty > 1$ and ${\cal
  I}_n(\tau) \neq \emptyset \ \forall n\in\N$ due to $({\cal P}1)$, Proposition
\ref{p.sna-existence} yields the existence of an SNA and an SNR. Further, we
have
\begin{eqnarray*}
  \Leb\left(\bigcup_{n=0}^{\infty} \bigcup_{k=-M_n-1}^{M_n+1} {\cal
      I}_n-k\omega \right)  & \leq & \sum_{n=0}^\infty 4N_n\eps_n \ \leq \ 
  \eps_0 N_0 +\sum_{n=1}^\infty \alpha^{-N_n/2p} \ .
\end{eqnarray*}
Again, the right side goes to $\Leb(\Lambda)$ as $\eps_0 \to 0$ and $\alpha \to
\infty$, such that (\ref{e.minimality-condition}) will be satisfied for small
$\eps_0$ and large $\alpha$. Consequently, we can apply
Proposition~\ref{p.minimality} to obtain 
 \ref{e.star}  for all $\tau\in\Lambda_\infty$. \smallskip

Finally, suppose that for some $\tau_0$ the symmetry condition
(\ref{e.symmetry-addendum}) holds. In this situation it follows by induction
that the critical regions ${\cal I}_n=I^1_n\cup I^2_n$ defined recursively by
(\ref{eq:defIn}) satisfy
\begin{equation} \label{e.interval-symmetry}
I_n^2\ =\ I^1_n+\halb \qquad \forall n\in\N \ .
\end{equation}
We show that in this case, for all sufficiently large $\alpha$, there exists a
sequence of integers $M_n\in[N_n,2N_n)$ such that the \ref{e.X'n} and \ref{e.Y'n}
hold for all $n\in\N$.  As before, Proposition~\ref{p.sna-existence} and
Proposition then imply \ref{e.star}, such that $\tau\in\Lambda_\infty$.

$({\cal X}')_0$ with $M_0=N_0=3$ holds for small $\eps_0$ due to the Diophantine
condition and $({\cal Y}')_0$ is void.  Suppose $M_0\ld M_n$ are chosen such
that \ref{e.X'n} and \ref{e.Y'n} hold, such that $\tau\in{\cal P}_n(M_0\ld
M_n)$. Then due to Lemma~\ref{l.Yn+1-Mexists} there exists
$M_{n+1}\in[N_{n+1},2N_{n+1})$ such that $({\cal Y}')_{n+1}$ holds.  Furthermore
$|I^\iota_{n+1}(\tau)| \leq \eps_{n+1}$ for $\iota=1,2$ due to $({\cal
  P}3)$. Now suppose that $({\cal X}')_{n+1}$ is not satisfied, such that
$I^\iota_{n+1}(\tau)\cap (I^\kappa_{n+1}(\tau) + n\omega) \neq \emptyset$ for
some $\iota,\kappa\in\{1,2\}$ and $|n|\leq 2K_{n+1}M_{n+1}$. Due to
(\ref{e.interval-symmetry}) this implies $d(2n\omega,0) \leq 2\eps_{n+1}$, which
contradicts the Diophantine condition~(\ref{eq:Diophantine}) when $\alpha$ is
large.  Consequently, when $\eps_0$ is sufficiently small and $\alpha$ is
sufficiently large conditions \ref{e.X'n} and \ref{e.Y'n} hold for all $n\in\N$
and we can apply Proposition \ref{p.minimality} to deduce that $f_{\tau_0}$
satisfies \ref{e.star}. Hence, $\tau_0$ can be included in
$\Lambda_\infty$. \qed

\section{The refined version of the twist parameter
  exclusion} \label{RefinedParameterExclusion}

The aim of this section is to prove Theorem~\ref{t.mr-quantitative2}. To that
end, we have to improve some of the estimates from the previous section by
taking into account the stronger assumptions on $\partial_\theta f_\theta$ in
(\ref{e.A8'}) and (\ref{eq:refinedbounddth}). As before, we can rely to some
extent on the respective results from \cite{jaeger:2009a}.

\subsection{Estimates on the critical sets and critical
  regions.}\label{PreviousRefinedEstimates}

Parts (i) and (ii) of Proposition~\ref{p.In-size} are replaced by the following
statements, which can again be taken from \cite{jaeger:2009a}.

\begin{prop}[Proposition 4.3 in \cite{jaeger:2009a}]
  \label{p.inductivelemma}
  Suppose (\ref{eq:Cinvariance})--(\ref{eq:s}), (\ref{eq:refinedbounddth}) and
  (\ref{e.S<Ad})--(\ref{e.s'<A}) hold, \ref{e.Xn} and \ref{e.Yn} are
  satisfied, $\alpha_\infty > 1$ and $M_0\geq d^{1/4}$. Then there exist
  a constant $d_1=d_1(\alpha_\infty)>0$ such that
  the following holds: If $d>d_1$ then \romanlist
\item ${\cal I}_n(\tau)$ and ${\cal I}_{n+1}(\tau)$ consist of exactly
  ${\cal N}$ connected components and each component of ${\cal I}_n(\tau)$
  contains exactly one component of ${\cal I}_{n+1}(\tau)$.
\item For all connected components of $I^\iota_{n+1}(\tau)$ of ${\cal
    I}_{n+1}(\tau)$ there holds $|I^\iota_{n+1}(\tau)| \leq 2\alpha_\infty^{-M_n}/s$.
\listend
\end{prop}
In contrast to this, the required version of Proposition~\ref{p.In-size}(iii)
has to take into account the fact that due to (\ref{e.A7'}) only two critical
regions exist. This assumption is not considered in \cite{jaeger:2009a}, such
that we cannot use the respective estimates there. Instead, we use the
following statement.
\begin{lem} \label{l.refined-slope-estimate} Suppose
  (\ref{eq:Cinvariance})--(\ref{eq:crossing}), (\ref{e.A7'})--(\ref{eq:refinedbounddth}) and
  (\ref{e.S<Ad})--(\ref{e.s'<A}) hold, \ref{e.Xn} and \ref{e.Yn} are
  satisfied, $\alpha_\infty > 1$ and $M_0\geq d^{1/4}$. Then there exist
  a constant $d_2=d_2(\alpha_\infty)>0$ such that
  the following holds. If $d>d_2$ then
  \begin{eqnarray}
    \label{e.theta-derivative1}
    s/2 & \leq & \partial_\theta(\varphi_{1,n}^\pm-\psi_{1,n}^\mp) \ \leq \ 2S 
    \qquad \ \ \textrm{ on } I^1_n(\tau)+\omega \ \textrm{ and } \\
    \label{e.theta-derivative2}
    -2S & \leq & \partial_\theta(\varphi_{2,n}^\pm-\psi_{2,n}^\mp) \ \leq \ -s/2 
   \qquad \textrm{ on } I^2_n(\tau)+\omega \ .
  \end{eqnarray}
\end{lem}
\myproof As in the proof of Lemma~\ref{l.lambda-dependence} and with the
notation introduced there, we have
\begin{equation}
  \partial_\theta \varphi^\pm_{\iota,n}(\theta,\tau) \ = \ \partial_\theta f_{\tau,\theta_{M_n-1}}(x_{M_{n}-1})
  + \underbrace{ \sum_{k=1}^{M_n-1} \partial_x f^{M_n-k}_{\tau,\theta_k}(x_k) \cdot 
   \partial_\theta f_{\tau,\theta_{k-1}}(x_{k-1})}_{=:q^\iota_1} \ .
\end{equation}
As $(\theta_{M_n},x_{M_n})=(\theta,\varphi^\pm_{\iota,n}(\theta,\tau))\in{\cal
I}_0+\omega$, we can use (\ref{eq:refinedbounddth}) together with
(\ref{e.finite-lyaps}) and (\ref{eq:bounddth}) to obtain
\begin{equation}
  |q_1^\iota| \ \leq \ \frac{s'+\alpha_\infty^{-M_0}S}{\alpha_\infty-1} \ .
\end{equation}
For $q^\iota_2:=\partial_\theta\psi^\mp_{\iota,n}(\theta,\tau)$ we obtain in a
similar way
\begin{equation}
  |q^\iota_2| \ \leq \ \frac{s'+\alpha_\infty^{-M_0}S}{\alpha_\infty-1} \ .
\end{equation}
Since $M_0\geq d^{1/4}$ we obtain that $|q^\iota_1|$ and $|q^\iota_2|$ are small
compared to $s$ and $S$ if $d$ is sufficiently large and $\frac{s'}{s}$ is
sufficiently small. As $\theta_{M_n-1}\in {\cal I}_0$, the statement follows
from (\ref{eq:bounddth}) and (\ref{e.A8'}).  \qed\medskip

In order to control the parameter dependence of the critical sets we replace
Lemma~\ref{l.lambda-dependence} by

 \begin{lem} \label{l.refined-tau-dependence} Suppose
   (\ref{eq:Cinvariance})--(\ref{eq:crossing}),
   (\ref{e.A7'})--(\ref{eq:refinedbounddth}) and
   (\ref{e.S<Ad})--(\ref{e.eps<1/Ad}) hold and \ref{e.Xn} and \ref{e.Yn} are
   satisfied. Further, let $d_2$ be chosen as in
   Lemma~\ref{l.refined-slope-estimate} and assume that $d>d_2$. Then
   \begin{eqnarray}
     \label{e.refined-tau-dependence}
     \partial_\tau I_{n+1}^1(\tau) \ \leq \ \frac{-\gamma}{2S} & \ , \ &
     \partial_\tau I^2_{n+1}(\tau) \ \geq \ \frac{\gamma}{2S} \\ \textrm{ and }
     \quad |\partial_\tau I^\iota_{n+1}(\tau)| & \leq &
     \frac{4L}{s(1-1/\alpha_\infty)} \quad (\iota=1,2) \ .
    \label{e.ref-tau-dep2}
   \end{eqnarray}
 \end{lem}
 \myproof Similar to the proof of Lemma~\ref{l.lambda-dependence} we
 let $I^\iota_{n+1}(\tau) = (a^\iota(\tau),b^\iota(\tau))$ and show that 
 \begin{equation}
   \label{e.ref-tau-dep3}
   -4L/s(1-1/\alpha_\infty) \ \leq \ \partial_\tau a^1(\tau) \ \leq -\gamma/2S \ .
 \end{equation}
 The required estimates on $\partial_\tau b^1(\tau),\ \partial_\tau
 a^2(\tau)$ and $\partial_\tau b^2(\tau)$ can then be treated in the
 same way. Note that (\ref{e.theta-derivative1}) in
 Lemma~\ref{l.refined-slope-estimate} implies that $f^{M_n}(A^1_n)$
 crosses $f^{-M_n}(B^1_n)$ upwards.

 We define $r^\iota_1$ and $r^\iota_2$ as in (\ref{e.rqiota}) and
 (\ref{e.psicontrol}). From (\ref{e.rqiota}) and (\ref{e.A8'}) we obtain that
\begin{equation}
  \partial_\tau\varphi^+_{\iota,n}(\theta,\tau) \ = \ \partial_\tau
f_{\tau,\theta_{M_n-1}}(x_{M_n-1})+r^\iota_1 \ \geq \ \partial_\tau
f_{\tau,\theta_{M_n-1}}(x_{M_n-1}) \ \geq \ \gamma \ .
\end{equation}
Note that $r^1_1\geq 0$ since all terms in the sum in (\ref{e.rqiota}) are
non-negative. Similarly $r^1_2\leq 0$, such that
\begin{equation}
  \partial_\tau(\varphi^+_{1,n}-\psi^-_{1,n})(a^1(\tau)+\omega,\tau) \ \geq \
  \gamma \ .
\end{equation}
Using (\ref{e.IFT-Formula}) and Lemma~\ref{l.refined-slope-estimate}
gives $\partial_\tau a^1(\tau) \leq -\gamma/2S$ as required. 

For the upper bound on $|\partial_\tau a^1(\tau)|$, note that using
(\ref{e.A8'}) and Lemma~\ref{p.lyaps} to estimate the sums defining $r^1_1$ and
$r^2_1$ in (\ref{e.rqiota}) and (\ref{e.psicontrol}) yields
\begin{equation}
  \partial_\tau(\varphi^+_{1,n}-\psi^-_{1,n})(a^1(\tau)+\omega,\tau) \
  \leq \ 2L/(1-1/\alpha_\infty) \ . 
\end{equation}
Together with Lemma~\ref{l.refined-slope-estimate}, this provides the
required bound $|\partial_\tau a^1(\tau)| \leq 4L/s(1-1/\alpha_\infty)$.
\qed \medskip

\begin{rem} \label{r.forcing-structure} The proof of
  Lemma~\ref{l.refined-tau-dependence} demonstrates well the restrictions which
  the need for controlling the relative speed of the critical intervals inflicts
  on the geormetry of the forcing. Considering the case of only two critical
  intervals with opposite sign of the slope of $\partial_\theta f_\theta$, as we
  do here, is not the only possibility to achive this. For instance, one could
  treat a multitude of critical intervals, as in Theorem~\ref{t.firstversion},
  by requiring that the twist $\partial_\tau f_\theta(x)$ almost vanishes
  outside of the critical regions (similar to the use of
  (\ref{eq:refinedbounddth}) in the proof of
  Lemma~\ref{l.refined-slope-estimate}). However, we see no way of treating more
  than two critical intervals if the twist is uniform as in
  (\ref{e.arnold}). The reason is that the lack of strong hyperbolicity does not
  allow to control the influence of the twist far from the critical region
  ${\cal I}_0(\tau)$ on the relative speed of the critical intervals. This could
  result in critical intervals moving at the same speed, in which case parameter
  exclusion would not work anymore. 
\end{rem}

\subsection{Proof of Theorem~\ref{t.mr-quantitative2}.} \label{RefinedProof}We
choose $t$ and $K_n=2^{n+t}/{\cal N}^2$ as in the proof of
Theorem~\ref{t.firstversion}, such that $\alpha_\infty\geq \alpha^{1/p}
>1$. Further, we suppose that $d$ satisfies $d\geq \max\{d_1,d_2\}$, where $d_1$
and $d_2$ are the constants from Proposition~\ref{p.inductivelemma} and
Lemma~\ref{l.refined-slope-estimate}. Then the mapping ${\cal I}_n:\Lambda\times
\N \to {\cal S}(\kreis)$ satisfies (\ref{e.P}), with $L$ replaced by
$4L/(1-1/\alpha_\infty)$ due to the weaker estimate on $|\partial_\tau
I^\iota_n(\tau)|$ in (\ref{e.refined-tau-dependence}) (compare
Remark~\ref{b.parameter-exclusion-assumption}).

Fix $k>2\nu$ and let $N_0$ be the first integer $\geq d^{1/k}$. As before, we
let $q=\max\{8,2\nu\}$ and define the sequences $N_n$ and $\eps_n$ recursively
by $N_{n+1}=\alpha^{N_n/qp}$ and $\eps_{n+1}=\frac{2}{s}\cdot
\alpha^{-N_n/p}$. Then using the dependencies (\ref{e.S<Ad})--(\ref{e.s'<A})
it is easy to check that all estimates on the quantities $v_{n},\ u_n$ and $V_n$
made in the proof of Theorem~\ref{t.firstversion} remain valid if the largeness
assumption on $\alpha$ used there is replaced by a largeness condition on $d$
that depends on the constants $\alpha,p,L,\gamma,A$ and $\delta$. Consequently,
the constant $m$ in Proposition~\ref{p.pex} satisfies
\begin{equation}
  m\ \geq \ \Leb(\Lambda) - u_0 - \sum_{n=0}^\infty \alpha^{-N_n/2p} \ .
\end{equation}
Furthermore, since $I^1_0(\tau)=I^2_0(\tau)+\halb$, $|I_0^\iota(\tau)|\leq
\eps_0\leq A/\sqrt{d}$ by (\ref{e.eps<1/Ad}) and $N_0 \leq d^{1/k}+1$, the
Diophantine condition implies that for sufficiently large $d$ condition $({\cal
X}')_0$ holds for all $\tau\in\Lambda$ (that is, ${\cal I}_0(\tau)$ is disjoint
from its first $2K_0N_0$ iterates).  This means that ${\cal P}(N_0)=\Lambda$ and
we can therefore apply Proposition~\ref{p.pex}(b), which yields a set
$\Lambda_\infty$ of measure
\[
\Leb(\Lambda_\infty) \ \geq \ \Leb(\Lambda) - \sum_{n=0}^\infty \alpha^{-N_n/2p}  \
\]
on which the slow-recurrence conditions \ref{e.X'n} and \ref{e.Y'n} hold for all
$n\in\N$. Consequently, for all $\tau\in\Lambda_\infty$ the existence of an SNA
and an SNR follows from Proposition~\ref{p.sna-existence}. Further, we have
\begin{eqnarray*}
  \Leb\left(\bigcup_{n=0}^{\infty} \bigcup_{k=-M_n-1}^{M_n+1} {\cal I}_n-k\omega
      \right) & \leq & \sum_{n=0}^\infty 2N_n\eps_n \ \leq \ \eps_0 N_0
      +\frac{2}{s}\cdot\sum_{n=1}^\infty \alpha^{-N_n/2p} \ .
\end{eqnarray*}
Due to (\ref{e.eps<1/Ad}) and the choice of $N_0$ in $[d^{1/k},d^{1/k}+1]$ the
sum on the right goes to zero as $d\to\infty$, and we can apply
Proposition~\ref{p.minimality} to obtain \ref{e.star}. Finally, the symmetry
statement is shown in the same way as in the proof of
Theorem~\ref{t.firstversion}.  \qed

\subsection{Proof of Corollaries~\ref{c.arnold} and \ref{t.arnold-half}.} 
\label{ArnoldCorollary} Recall that we consider the parameter family
\[
f_{a,b,\tau}(\theta,x) \ = \ \left(\theta+\omega,x+\tau+\frac{a}{2\pi}\sin(2\pi x) +
g_b(\theta)\right) \ 
\]
with $g_b(\theta)=\arctan(b\sin(2\pi \theta))/\pi$. We suppose that
$\omega$ is Diophantine with constants $\sigma,\nu$, such that
(\ref{eq:Diophantine}) holds. Let $h_a(x) = x+\frac{a}{2\pi}\sin(2\pi x)$.
Then there exist constants $0<e<c<\halb$ $\alpha>1,\ p\in\N$ and
$0<t<\halb-c$ such that there holds
\begin{eqnarray}
  h_a\left([e,-e]\right) & \ssq & (c+t,-c-t) \ ,\\
  h_a'(x) & < & \alpha^{-2/p} \qquad \forall x\in C:=[c,-c] \ , \\
  h_a'(x) & > & \alpha^{2/p} \qquad \forall x\in E:= [-e,e] \quad \textrm{and} \\
  \alpha^{-p} & < & h_a'(x) \ < \ \alpha^{p} \qquad \forall x\in\kreis \ . 
\end{eqnarray}
Let $\Lambda := \left(\halb-\frac{t}{2},\halb+\frac{t}{2}\right)$ and
$\gamma_0:=\viertel\tan(\pi(\halb-\frac{t}{2}))$. Further, define
$I^1_0:=\left(-\gamma_0/b,\gamma_0/b\right)$ and $I^2_0:=I^1_0+\halb$. Then for
$\theta\notin {\cal I}_0=I^1_0\cup I^2_0$ and $\tau\in\Lambda$ we obtain
\[\textstyle
d(g_b(\theta)+\tau,0) \ \leq \ \frac{t}{2} +d(g_b(\theta),\halb) \ < t \ .
\]
Consequently $f_{a,b,\tau}$ satisfies
(\ref{eq:Cinvariance})--(\ref{eq:bounds3}) for all $\tau\in\Lambda$. Further, we have 
\begin{eqnarray*}
  \partial_\theta f_\theta(x) & = & \partial_\theta g_b(\theta) \ = \
  \frac{2}{1+b^2\sin^2(2\pi\theta)} \cdot b\cos(2\pi\theta) \ \leq \ 2b
  \qquad \forall (\theta,x)\in\torus \quad \textrm{and}\\
  |\partial_\theta f_\theta(x)| & \geq & \frac{b}{1+(2\pi\gamma_0)^2} \qquad 
  \forall(\theta,x) \in {\cal I}_0\times \kreis \ .
\end{eqnarray*}
This allows to see that (\ref{eq:bounddth}), (\ref{eq:crossing}), (\ref{e.A7'})
and (\ref{e.A8'}) are satisfied with $S=2b$ and
$s=\frac{b}{1+(2\pi\gamma_0)^2}$. If we let
$I^{1'}_0=[-\gamma_0/\sqrt{b},\gamma/\sqrt{b}]$, $I^2_0=I^{1'}_0+\halb$ and
${\cal I}_0'=I^{1'}_0\cup I^{2'}_0$ then
\[
|\partial_\theta f_\theta(x)| \ \leq \ 2\gamma_0^{-2}
\qquad \forall (\theta,x)\in(\kreis \smin {\cal I}_0')\times\kreis \ ,
\]
such that (\ref{eq:refinedbounddth}) holds with
$s'=2\gamma_0^{-2}$. Finally, since $\partial_\tau f_\theta(x)
= 1$, we may choose $\gamma=\halb$ and $L=2$ in (\ref{e.A10'}).  Altogether,
this implies that all assumptions of Theorem~\ref{t.mr-quantitative2} are
satisfied for a suitable constant $A$ and $d=b$. The conclusions of the corollaries follow. \qed


\begin{thebibliography}{10}
\small

\bibitem{benedicks/carleson:1991}
M.~Benedicks and L.~Carleson.
\newblock The dynamics of the {H\'enon} map.
\newblock {\em Ann. Math. (2)}, 133(1):73--169, 1991.

\bibitem{milnor:1985}
J.~Milnor.
\newblock On the concept of attractor.
\newblock {\em Commun. Math. Phys.}, 99:177--195, 1985.

\bibitem{grebogi/ott/pelikan/yorke:1984}
C.~Grebogi, E.~Ott, S.~Pelikan, and J.A. Yorke.
\newblock Strange attractors that are not chaotic.
\newblock {\em Physica D}, 13:261--268, 1984.

\bibitem{keller:1996}
G.~Keller.
\newblock A note on strange nonchaotic attractors.
\newblock {\em Fundam.\ Math.}, 151(2):139--148, 1996.

\bibitem{romeiras/etal:1987}
F.J. Romeiras, A.~Bondeson, E.~Ott, T.M. Antonsen~Jr., and C.~Grebogi.
\newblock Quasiperiodically forced dynamical systems with strange nonchaotic
  attractors.
\newblock {\em Physica D}, 26:277--294, 1987.

\bibitem{ding/grebogi/ott:1989}
M.~Ding, C.~Grebogi, and E.~Ott.
\newblock Evolution of attractors in quasiperiodically forced systems: {F}rom
  quasiperiodic to strange nonchaotic to chaotic.
\newblock {\em Phys. Rev. A}, 39(5):2593--2598, 1989.

\bibitem{feudel/kurths/pikovsky:1995}
U.~Feudel, J.~Kurths, and A.~Pikovsky.
\newblock Strange nonchaotic attractor in a quasiperiodically forced circle
  map.
\newblock {\em Physica D}, 88:176--186, 1995.

\bibitem{millionscikov:1969}
V.M. Million{\u{s}\u{c}}ikov.
\newblock Proof of the existence of irregular systems of linear differential
  equations with quasi periodic coefficients.
\newblock {\em Differ.~Uravn.}, 5(11):1979--1983, 1969.

\bibitem{vinograd:1975}
R.E Vinograd.
\newblock A problem suggested by {N.R.~Erugin}.
\newblock {\em Differ.~Uravn.}, 11(4):632--638, 1975.

\bibitem{herman:1983}
M.~Herman.
\newblock Une {m\'{e}thode} pour minorer les exposants de {{L}yapunov} et
  quelques exemples montrant le {caract\`{e}re} local d'un {th\'{e}or\`{e}me}
  {d'Arnold} et de {Moser} sur le tore de dimension 2.
\newblock {\em Comment.\ Math.\ Helv.}, 58:453--502, 1983.

\bibitem{avila/krikorian:2004}
A.~Avila and R.~Krikorian.
\newblock Reducibility or non-uniform hyperbolicity for quasiperiodic
  {S}chr{\"o}dinger cocycles.
\newblock {\em Ann. Math. (2)}, 164:911--940, 2006.

\bibitem{haro/puig:2006}
A.~Haro and J.~Puig.
\newblock Strange non-chaotic attractors in {H}arper maps.
\newblock {\em Chaos}, 16, 2006.

\bibitem{puig:2004}
J.~Puig.
\newblock Cantor spectrum for the almost {M}athieu operator.
\newblock {\em Comm. Math. Phys.}, 244(2):297--309, 2004.

\bibitem{avila/jitomirskaya:2005}
A.~Avila and S.~Jitomirskaya.
\newblock The {T}en {M}artini {P}roblem.
\newblock {\em Ann. Math. (2)}, 170(1):303--342, 2009.

\bibitem{avila/jitomirskaya:2010}
A.~Avila and S.~Jitomirskaya.
\newblock Almost localization and almost deducibility.
\newblock {\em J. Eur. Math. Soc.}, 12(1):93--131, 2010.

\bibitem{avila:2010}
A.~Avila.
\newblock Global theory of one-frequency {Schr\"o}dinger operators {I and II}.
\newblock Preprints 2010.

\bibitem{young:1997}
L.-S. Young.
\newblock Lyapunov exponents for some quasi-periodic cocycles.
\newblock {\em Ergodic Theory Dyn. Syst.}, 17:483--504, 1997.

\bibitem{bjerkloev:2005a}
K.~Bjerkl{\"o}v.
\newblock Positive {L}yapunov exponent and minimality for a class of
  one-dimensional quasi-periodic {S}chr{\"o}dinger equations.
\newblock {\em Ergodic Theory Dyn. Syst.}, 25:1015--1045, 2005.

\bibitem{jaeger:2009a}
T.~J{\"a}ger.
\newblock Strange non-chaotic attractors in quasiperiodically forced circle
  maps.
\newblock {\em Comm. Math. Phys.}, 289(1):253--289, 2009.

\bibitem{jaeger:2006a}
T.~J{\"a}ger.
\newblock The creation of strange non-chaotic attractors in non-smooth
  saddle-node bifurcations.
\newblock {\em Mem.\ Am.\ Math.\ Soc.}, 945:1--106, 2009.

\bibitem{bjerkloev/jaeger:2009}
K.~Bjerkl{\"o}v and T.~J{\"a}ger.
\newblock Rotation numbers for quasiperiodically forced circle maps --
  {M}ode-locking vs strict monotonicity.
\newblock {\em J.\ Am.\ Math.\ Soc.}, 22(2):353--362, 2009.

\bibitem{bellissard/simon:1982}
J.~B{\a'e}llissard and B.~Simon.
\newblock Cantor spectrum for the almost {Mathieu} equation.
\newblock {\em J. Funct. Anal.}, 48(3):408--419, 1982.

\bibitem{stark/feudel/glendinning/pikovsky:2002}
J.~Stark, U.~Feudel, P.~Glendinning, and A.~Pikovsky.
\newblock Rotation numbers for quasi-periodically forced monotone circle maps.
\newblock {\em Dyn.\ Syst.}, 17(1):1--28, 2002.

\end{thebibliography}

\end{document}